\documentclass[12pt,english]{article}
\usepackage{avant}
\usepackage[T1]{fontenc}
\usepackage[utf8]{inputenc}
\usepackage{geometry}
\usepackage{babel}
\usepackage{amsmath}
\usepackage{amssymb}
\usepackage{stmaryrd}
\usepackage[authoryear]{natbib}
\usepackage{algorithm}
\usepackage{graphicx}

\usepackage{float}
\usepackage{booktabs} 
\usepackage{setspace}
\usepackage{subcaption}
\usepackage{multirow}
\usepackage{multicol}
\usepackage{threeparttable}
\usepackage[unicode=true,
bookmarks=false,
breaklinks=false,pdfborder={0 0 1},backref=false,colorlinks=false]
{hyperref}
\usepackage{breakurl}

\makeatletter
\@ifundefined{date}{}{\date{}}

\usepackage{array}
\newcolumntype{L}[1]{>{\raggedright\let\newline\\\arraybackslash\hspace{0pt}}m{#1}}
\newcolumntype{C}[1]{>{\centering\let\newline\\\arraybackslash\hspace{0pt}}m{#1}}
\newcolumntype{R}[1]{>{\raggedleft\let\newline\\\arraybackslash\hspace{0pt}}m{#1}}
\usepackage{float}
\usepackage{multirow}
\usepackage{breakurl}

\providecommand{\tabularnewline}{\\}
\floatstyle{ruled}

\usepackage{xr}
\makeatletter
\newcommand*{\addFileDependency}[1]{
	\typeout{(#1)}
	\@addtofilelist{#1}
	\IfFileExists{#1}{}{\typeout{No file #1.}}
}
\makeatother

\providecommand{\tabularnewline}{\\}
\floatstyle{ruled}

\usepackage{amsthm}\usepackage{bm}\usepackage{rotating}\usepackage{array}\usepackage{booktabs}\usepackage{color}
\usepackage{multirow}\usepackage{enumerate}\usepackage{caption}
\usepackage{float}
\usepackage{multirow}
\providecommand{\tabularnewline}{\\}
\floatstyle{ruled}
\newfloat{algorithm}{tbp}{loa}
\providecommand{\algorithmname}{Algorithm}
\floatname{algorithm}{\protect\algorithmname}

\usepackage{bm}\usepackage{rotating}\usepackage{array}\usepackage{booktabs}\usepackage{color}\usepackage{algorithmic}
\usepackage{algorithm}\usepackage{multirow}\usepackage{enumerate}\usepackage{caption}\usepackage{subcaption}
\floatstyle{ruled}
\newfloat{algorithm}{tbp}{loa}
\providecommand{\algorithmname}{Algorithm}
\floatname{algorithm}{\protect\algorithmname}

\topmargin by -\headheight \advance
\topmargin by -\headsep     

\evensidemargin
\oddsidemargin
\providecommand{\tabularnewline}{\\}

\RequirePackage[displaymath]{lineno}

\usepackage{amsfonts}\usepackage{latexsym}

\usepackage{mathrsfs}\usepackage{amsthm}\usepackage{latexsym}\usepackage{booktabs}

\usepackage{times}

\DeclareMathOperator*{\argmax}{argmax}

\newtheorem{lemma}{{\bf Lemma}}

\newtheorem{theorem}{{\bf Theorem}}

\makeatother
\begin{document}
	
	\renewcommand{\thepage}{\arabic{page}}  
	\renewcommand{\thesection}{\arabic{section}}   
	\renewcommand{\thetable}{\arabic{table}}   
	\renewcommand{\thefigure}{\arabic{figure}}
	\renewcommand{\theequation}{\arabic{equation}}
	\renewcommand{\thetheorem}{\arabic{theorem}}
	\renewcommand{\thelemma}{\arabic{lemma}}
	\renewcommand{\thealgorithm}{\arabic{algorithm}}
	{\setlength{\baselineskip}{1.5\baselineskip} 
\global\long\def\mba{{a}}
\global\long\def\mbA{{A}}
\global\long\def\mbb{{b}}
\global\long\def\mbB{{B}}
\global\long\def\mbc{{c}}
\global\long\def\mbC{{C}}
\global\long\def\mbd{{d}}
\global\long\def\mbD{{D}}
\global\long\def\mbe{{e}}
\global\long\def\mbE{{E}}
\global\long\def\mbf{{f}}
\global\long\def\mbF{{F}}
\global\long\def\mbg{{g}}
\global\long\def\mbG{{G}}
\global\long\def\mbh{{h}}
\global\long\def\mbH{{H}}
\global\long\def\mbi{{i}}
\global\long\def\mbI{{I}}
\global\long\def\mbj{{j}}
\global\long\def\mbJ{{J}}
\global\long\def\mbk{{k}}
\global\long\def\mbK{{K}}
\global\long\def\mbl{{l}}
\global\long\def\mbL{{L}}
\global\long\def\mbm{{m}}
\global\long\def\mbM{{M}}
\global\long\def\mbn{{n}}
\global\long\def\mbN{{N}}
\global\long\def\mbo{{o}}
\global\long\def\mbO{{O}}
\global\long\def\mbp{{p}}
\global\long\def\mbP{{P}}
\global\long\def\mbq{{q}}
\global\long\def\mbQ{{Q}}
\global\long\def\mbr{{r}}
\global\long\def\mbR{{R}}
\global\long\def\mbs{{s}}
\global\long\def\mbS{{S}}
\global\long\def\mbt{{t}}
\global\long\def\mbT{{T}}
\global\long\def\mbu{{u}}
\global\long\def\mbU{{U}}
\global\long\def\mbv{{v}}
\global\long\def\mbV{{V}}
\global\long\def\mbw{{w}}
\global\long\def\mbW{{W}}
\global\long\def\mbx{{x}}
\global\long\def\mbX{{X}}
\global\long\def\mby{{y}}
\global\long\def\mbY{{Y}}
\global\long\def\mbz{{z}}
\global\long\def\mbZ{{Z}}

\global\long\def\hatmba{\widehat{{a}}}
\global\long\def\hatmbA{\widehat{{A}}}
\global\long\def\hatmbb{\widehat{{b}}}
\global\long\def\hatmbB{\widehat{{B}}}
\global\long\def\hatmbc{\widehat{{c}}}
\global\long\def\hatmbC{\widehat{{C}}}
\global\long\def\hatmbd{\widehat{{d}}}
\global\long\def\hatmbD{\widehat{{D}}}
\global\long\def\hatmbe{\widehat{{e}}}
\global\long\def\hatmbE{\widehat{{E}}}
\global\long\def\hatmbf{\widehat{{f}}}
\global\long\def\hatmbF{\widehat{{F}}}
\global\long\def\hatmbg{\widehat{{g}}}
\global\long\def\hatmbG{\widehat{{G}}}
\global\long\def\hatmbh{\widehat{{h}}}
\global\long\def\hatmbH{\widehat{{H}}}
\global\long\def\hatmbi{\widehat{{i}}}
\global\long\def\hatmbI{\widehat{{I}}}
\global\long\def\hatmbj{\widehat{{j}}}
\global\long\def\hatmbJ{\widehat{{J}}}
\global\long\def\hatmbk{\widehat{{k}}}
\global\long\def\hatmbK{\widehat{{K}}}
\global\long\def\hatmbl{\widehat{{l}}}
\global\long\def\hatmbL{\widehat{{L}}}
\global\long\def\hatmbm{\widehat{{m}}}
\global\long\def\hatmbM{\widehat{{M}}}
\global\long\def\hatmbn{\widehat{{n}}}
\global\long\def\hatmbN{\widehat{{N}}}
\global\long\def\hatmbo{\widehat{{o}}}
\global\long\def\hatmbO{\widehat{{O}}}
\global\long\def\hatmbp{\widehat{{p}}}
\global\long\def\hatmbP{\widehat{{P}}}
\global\long\def\hatmbq{\widehat{{q}}}
\global\long\def\hatmbQ{\widehat{{Q}}}
\global\long\def\hatmbr{\widehat{{r}}}
\global\long\def\hatmbR{\widehat{{R}}}
\global\long\def\hatmbs{\widehat{{s}}}
\global\long\def\hatmbS{\widehat{{S}}}
\global\long\def\hatmbt{\widehat{{t}}}
\global\long\def\hatmbT{\widehat{{T}}}
\global\long\def\hatmbu{\widehat{{u}}}
\global\long\def\hatmbU{\widehat{{U}}}
\global\long\def\hatmbv{\widehat{{v}}}
\global\long\def\hatmbV{\widehat{{V}}}
\global\long\def\hatmbw{\widehat{{w}}}
\global\long\def\hatmbW{\widehat{{W}}}
\global\long\def\hatmbx{\widehat{{x}}}
\global\long\def\hatmbX{\widehat{{X}}}
\global\long\def\hatmby{\widehat{{y}}}
\global\long\def\hatmbY{\widehat{{Y}}}
\global\long\def\hatmbz{\widehat{{z}}}
\global\long\def\hatmbZ{\widehat{{Z}}}

\global\long\def\bolalpha{{\alpha}}
\global\long\def\bolbeta{{\beta}}
\global\long\def\bolgamma{{\gamma}}
\global\long\def\boldelta{{\delta}}
\global\long\def\bolepsilon{{\epsilon}}
\global\long\def\bolvarepsilon{{\varepsilon}}
\global\long\def\bolzeta{{\zeta}}
\global\long\def\boleta{{\eta}}
\global\long\def\boltheta{{\theta}}
\global\long\def\bolkappa{{\kappa}}
\global\long\def\bollambda{{\lambda}}
\global\long\def\bolmu{{\mu}}
\global\long\def\bolnu{{\nu}}
\global\long\def\bolxi{{\xi}}
\global\long\def\bolXi{{\Xi}}
\global\long\def\bolpi{{\pi}}
\global\long\def\bolrho{{\rho}}
\global\long\def\bolsigma{{\sigma}}
\global\long\def\boltau{{\tau}}
\global\long\def\bolphi{{\phi}}
\global\long\def\bolchi{{\chi}}
\global\long\def\bolpsi{{\psi}}
\global\long\def\bolomega{{\omega}}
\global\long\def\bolGamma{{\Gamma}}
\global\long\def\bolDelta{{\Delta}}
\global\long\def\bolTheta{{\Theta}}
\global\long\def\bolLambda{{\Lambda}}
\global\long\def\bolPi{{\Pi}}
\global\long\def\bolSigma{{\Sigma}}
\global\long\def\bolPhi{{\Phi}}
\global\long\def\bolPsi{{\Psi}}
\global\long\def\bolOmega{{\Omega}}

\global\long\def\hatbolalpha{\widehat{{\alpha}}}
\global\long\def\hatbolbeta{\widehat{{\beta}}}
\global\long\def\hatbolgamma{\widehat{{\gamma}}}
\global\long\def\hatboldelta{\widehat{{\delta}}}
\global\long\def\hatbolepsilon{\widehat{{\epsilon}}}
\global\long\def\hatbolzeta{\widehat{{\zeta}}}
\global\long\def\hatboleta{\widehat{{\eta}}}
\global\long\def\hatboltheta{\widehat{{\theta}}}
\global\long\def\hatbolkappa{\widehat{{\kappa}}}
\global\long\def\hatbollambda{\widehat{{\lambda}}}
\global\long\def\hatbolmu{\widehat{{\mu}}}
\global\long\def\hatbolnu{\widehat{{\nu}}}
\global\long\def\hatbolxi{\widehat{{\xi}}}
\global\long\def\hatbolXi{\widehat{{\Xi}}}
\global\long\def\hatbolpi{\widehat{{\pi}}}
\global\long\def\hatbolrho{\widehat{{\rho}}}
\global\long\def\hatbolsigma{\widehat{{\sigma}}}
\global\long\def\hatboltau{\widehat{{\tau}}}
\global\long\def\hatbolphi{\widehat{{\phi}}}
\global\long\def\hatbolchi{\widehat{{\chi}}}
\global\long\def\hatbolpsi{\widehat{{\psi}}}
\global\long\def\hatbolomega{\widehat{{\omega}}}
\global\long\def\hatbolGamma{\widehat{{\Gamma}}}
\global\long\def\hatbolDelta{\widehat{{\Delta}}}
\global\long\def\hatbolTheta{\widehat{{\Theta}}}
\global\long\def\hatbolLambda{\widehat{{\Lambda}}}
\global\long\def\hatbolPi{\widehat{{\Pi}}}
\global\long\def\hatbolSigma{\widehat{{\Sigma}}}
\global\long\def\hatbolPhi{\widehat{{\Phi}}}
\global\long\def\hatbolPsi{\widehat{{\Psi}}}
\global\long\def\hatbolOmega{\widehat{{\Omega}}}

\global\long\def\tilbolalpha{\widetilde{{\alpha}}}
\global\long\def\tilbolbeta{\widetilde{{\beta}}}
\global\long\def\tilbolgamma{\widetilde{{\gamma}}}
\global\long\def\tilboldelta{\widetilde{{\delta}}}
\global\long\def\tilbolepsilon{\widetilde{{\epsilon}}}
\global\long\def\tilbolzeta{\widetilde{{\zeta}}}
\global\long\def\tilboleta{\widetilde{{\eta}}}
\global\long\def\tilboltheta{\widetilde{{\theta}}}
\global\long\def\tilbolkappa{\widetilde{{\kappa}}}
\global\long\def\tilbollambda{\widetilde{{\lambda}}}
\global\long\def\tilbolmu{\widetilde{{\mu}}}
\global\long\def\tilbolnu{\widetilde{{\nu}}}
\global\long\def\tilbolxi{\widetilde{{\xi}}}
\global\long\def\tilbolpi{\widetilde{{\pi}}}
\global\long\def\tilbolrho{\widetilde{{\rho}}}
\global\long\def\tilbolsigma{\widetilde{{\sigma}}}
\global\long\def\tilboltau{\widetilde{{\tau}}}
\global\long\def\tilbolphi{\widetilde{{\phi}}}
\global\long\def\tilbolchi{\widetilde{{\chi}}}
\global\long\def\tilbolpsi{\widetilde{{\psi}}}
\global\long\def\tilbolomega{\widetilde{{\omega}}}
\global\long\def\tilbolGamma{\widetilde{{\Gamma}}}
\global\long\def\tilbolDelta{\widetilde{{\Delta}}}
\global\long\def\tilbolTheta{\widetilde{{\Theta}}}
\global\long\def\tilbolLambda{\widetilde{{\Lambda}}}
\global\long\def\tilbolPi{\widetilde{{\Pi}}}
\global\long\def\tilbolSigma{\widetilde{{\Sigma}}}
\global\long\def\tilbolPhi{\widetilde{{\Phi}}}
\global\long\def\tilbolPsi{\widetilde{{\Psi}}}
\global\long\def\tilbolOmega{\widetilde{{\Omega}}}

\global\long\def\barbolmu{\overline{\bolmu}}
\global\long\def\barmbX{\overline{\mbX}}

\global\long\def\td{\textsl{t}}

\global\long\def\mbbR{\mathbb{R}}
\global\long\def\mbbP{\mathbb{P}}
\global\long\def\mbbQ{\mathbb{Q}}
\global\long\def\mbbS{\mathbb{S}}
\global\long\def\mbbH{\mathbb{H}}
\global\long\def\mbbX{\mathbb{X}}
\global\long\def\mbbY{\mathbb{Y}}
\global\long\def\mbbW{\mathbb{W}}
\global\long\def\mbbZ{\mathbb{Z}}
\global\long\def\spc{\mathcal{S}}

\global\long\def\calA{\mathcal{A}}
\global\long\def\calB{\mathcal{B}}
\global\long\def\calC{\mathcal{C}}
\global\long\def\calD{\mathcal{D}}
\global\long\def\calE{\mathcal{E}}
\global\long\def\calF{\mathcal{F}}
\global\long\def\calG{\mathcal{G}}
\global\long\def\calH{\mathcal{H}}
\global\long\def\calI{\mathcal{I}}
\global\long\def\calJ{\mathcal{J}}
\global\long\def\calK{\mathcal{K}}
\global\long\def\calL{\mathcal{L}}
\global\long\def\calM{\mathcal{M}}
\global\long\def\calN{\mathcal{N}}
\global\long\def\calO{\mathcal{O}}
\global\long\def\calP{\mathcal{P}}
\global\long\def\calQ{\mathcal{Q}}
\global\long\def\calR{\mathcal{R}}
\global\long\def\calS{\mathcal{S}}
\global\long\def\calT{\mathcal{T}}
\global\long\def\calU{\mathcal{U}}
\global\long\def\calV{\mathcal{V}}
\global\long\def\calW{\mathcal{W}}

\global\long\def\mbell{{\ell}}
\global\long\def\bolell{{\ell}}
\global\long\def\mbzero{\mathbf{0}}

\global\long\def\bolPhio{{\Phi}_{0}}
\global\long\def\bolOmegao{{\Omega}_{0}}

\global\long\def\bolSigmaX{\bolSigma_{\mbX}}
\global\long\def\bolSigmaY{\bolSigma_{\mbY}}
\global\long\def\bolSigmaXY{{\Sigma}_{\mbX\mbY}}
\global\long\def\mbSX{\mathbf{S}_{\mbX}}
\global\long\def\mbSY{\mathbf{S}_{\mbY}}
\global\long\def\mbSXY{\mathbf{S}_{\mbX\mbY}}
\global\long\def\mbSYX{\mathbf{S}_{\mbY\mbX}}
\global\long\def\mbRYX{\mathbf{S}_{\mbY|\mbX}}
\global\long\def\mbRXY{\mathbf{S}_{\mbX|\mbY}}
\global\long\def\mbSc{\mbS_{\mbC}}
\global\long\def\mbSd{\mbS_{\mbD}}

\global\long\def\sumn{\sum_{i=1}^{n}}

\global\long\def\E{\mathrm{E}}

\global\long\def\APT{\mathrm{APT}}
\global\long\def\OST{\mathrm{OST}}
\global\long\def\OS{\mathrm{OS}}
\global\long\def\AR{\mathrm{AR}}
\global\long\def\APN{\mathrm{APN}}
\global\long\def\APL{\mathrm{APL}}
\global\long\def\REE{\mathrm{REE}}
\global\long\def\TPR{\mathrm{TPR}}
\global\long\def\FPR{\mathrm{FPR}}
\global\long\def\AGPT{\mathrm{AGPT}}

\global\long\def\TN{\mathrm{TN}}
\global\long\def\HOST{\mathrm{HOST}}
\global\long\def\TT{\mathrm{TT}}
\global\long\def\MVT{\mathrm{MVT}}

\global\long\def\L{\mathrm{L}}
\global\long\def\F{\mathrm{F}}
\global\long\def\J{\mathrm{J}}
\global\long\def\H{\mathrm{H}}
\global\long\def\G{\mathrm{G}}
\global\long\def\Cov{\mathrm{cov}}
\global\long\def\cov{\mathrm{cov}}
\global\long\def\Corr{\mathrm{corr}}
\global\long\def\Var{\mathrm{var}}
\global\long\def\dimension{\mathrm{dim}}
\global\long\def\spn{\mathrm{span}}
\global\long\def\vech{\mathrm{vech}}
\global\long\def\vecc{\mathrm{vec}}
\global\long\def\Prob{\mathrm{Pr}}
\global\long\def\Env{\mathrm{env}}
\global\long\def\tr{\mathrm{tr}}
\global\long\def\dg{\mathrm{diag}}
\global\long\def\asyVar{\mathrm{avar}}
\global\long\def\MSE{\mathrm{MSE}}
\global\long\def\OLS{\mathrm{OLS}}
\global\long\def\LDA{\mathrm{LDA}}
\global\long\def\QDA{\mathrm{QDA}}
\global\long\def\logg{\mathrm{log}}
\global\long\def\exp{\mathrm{exp}}
\global\long\def\floor{\mathrm{floor}}
\global\long\def\tr{\mathrm{tr}}
\global\long\def\spann{\mathrm{span}}
\global\long\def\TEC{\mathrm{TEC}}
\global\long\def\EC{\mathrm{EC}}
\global\long\def\N{\mathrm{N}}

\global\long\def\sgn{\mathrm{sgn}}
\global\long\def\CS{\calS_{Y\mid\mbX}}
\global\long\def\covenv{\calE_{\bolDelta_{0}}(\bolbeta)}
\global\long\def\mt{\mathrm{t}}
\global\long\def\mW{\mathrm{W}}
\global\long\def\Wis{\mathrm{Wishart}}
		
		\global\long\def\CS{\calS_{Y\mid\mbX}}
		\global\long\def\covenv{\calE_{\bolDelta_{0}}(\bolbeta)}
		\title{{Statistical analysis for a penalized EM algorithm in high-dimensional mixture linear regression model}}
		\author{Ning Wang, Xin Zhang and Qing Mai\thanks{Ning Wang ({ningwangbnu@bnu.edu.cn}) is Assistant Professor, Center of Statistics and Data Science, Beijing Normal University, Zhuhai, 519807, China; Xin Zhang (xzhang8@fsu.edu) is Associate Professor,  Department of Statistics, Florida State University, Tallahassee, 32312, Florida, USA; Qing Mai (qmai@fsu.edu) is Associate Professor,  Department of Statistics, Florida State University, Tallahassee, 32312, Florida, USA.}}
		\date{Beijing Normal University and Florida State University}
		\maketitle
		
		\begin{abstract}
			The expectation-maximization (EM) algorithm and its variants are widely used in statistics. In high-dimensional mixture linear regression, the model is assumed to be a finite mixture of linear regression and the number of predictors is much larger than the sample size. The standard EM algorithm, which attempts to find the maximum likelihood estimator, becomes infeasible for such model. We devise a group lasso penalized EM algorithm and study its statistical properties. Existing theoretical results of regularized EM algorithms often rely on dividing the sample into many independent batches and employing a fresh batch of sample in each iteration of the algorithm. Our algorithm and theoretical analysis do not require sample-splitting, and can be extended to multivariate response cases. The proposed methods also have encouraging performances in numerical studies.\\
			
		\end{abstract}
		
		\section{Introduction}
	Consider a univariate response $Y\in\mbbR$ and a $p$-dimensional predictor $X\in\mbbR^p$. 
	The mixture linear regression model assumes that
	\begin{equation}\label{mlr}
		Y=\beta_k^T X+\epsilon,\ \mathrm{for}\ k=1,\cdots,K,\ \mathrm{with}\ \mathrm{probability}\ \omega_k>0,
	\end{equation}
	where $K\geq 2$ is the number of mixtures, $\sum_{k=1}^K\omega_k=1$, $\epsilon\sim N(0,\sigma^2)$, $\sigma^2>0$, is independent of $X$, and $\beta_k$ is the $p$-dimensional regression coefficient vector that characterizes the linear relationship between $Y$ and $X$ in the $k$-th mixture. By introducing a latent variable $W\in \{1,\cdots,K\}$, independent of $X$, model \eqref{mlr} is equivalent to
	\begin{equation}\label{mlr1}
		\mbbP(W=k)=\omega_k,\quad Y\mid (X,W=k)\sim N(\beta_k^T X,\sigma^2).
	\end{equation}
	
	
	We consider the high-dimensional joint estimation of all the $\bolbeta_k$'s and provide a general estimation procedure with strong theoretical guarantees. The latent mixtures, indicated by the latent variable $W$ in \eqref{mlr1}, make the estimation problem much more challenging than the linear regression, especially in high dimensions. We assume that only a subset of predictors, indexed by $S\subset\{1,\cdots,p\}$, is relevant to the regression. Therefore, $(\beta_k)_{\calS^c}=0$ for all $k$, where $S^c$ is the complement of $S$. The group lasso penalty \citep{yuan2006model} is naturally applied to those $p$-dimensional $\beta_k$ vectors in the maximization steps of our regularized EM algorithm to select relevant variables across all the mixtures. 
	
	
	The mixture linear model and finite mixture models in general, are widely used to account for heterogeneity in data analysis \citep[e.g.,][]{turner2000estimating,mclachlan2004finite,mclachlan2019finite}. When the number of predictors is not large, the latent mixtures and the model parameters can be estimated using the expectation-maximization \citep[EM,][]{dempster1977maximum} algorithm. The EM algorithm is the dominant solution for finding the maximum likelihood estimators of mixture regression models like \eqref{mlr}. When $\epsilon$ is not normally distributed, the EM algorithm has been modified and extended to robust fitting of mixture linear models under the  t- and the Laplace distributions in \citet{yao2014robust} and \citet{song2014robust}, respectively.  \cite{flexmix} provides computational and implementation details of mixture regression models. 
	
	Although the EM algorithm has been extensively used in mixture regression models, it is challenging to establish a rigorous theoretical characterization of the finite-sample estimates in the iterative algorithm. Some groundbreaking progress has been made in recent years. \cite{balakrishnan2017statistical} laid theoretical foundations for quantifying the EM updates' convergence within statistical precision of a global optimum. For mixture linear model, strong theoretical guarantees of the EM algorithm are often established by focusing on the model with two equal mixtures, $K=2$ and $\omega_1=\omega_2=1/2$, and symmetric regression coefficient vectors, $\beta_2=-\beta_1$, \citep[e.g.,][and references therein]{kwon2020minimax}. Then the EM algorithm is simplified substantially because $\omega_1=\omega_2=1/2$ does not require estimation and, more importantly, the model parameters are reduced to a single vector $\beta\equiv\beta_1=-\beta_2\in\mbbR^p$. Many theoretical results are established for the ``sample-splitting'' EM algorithm, which divides the full data into $T$ equal batches and uses a new batch of samples in each iteration. Without sample splitting, theoretical analysis becomes much more challenging. This is because the function to be maximized in each EM iteration, namely the Q-function, involves both the random samples and the current parameter estimates, which are made independent by sample-splitting. See \cite{balakrishnan2017statistical} and \cite{klusowski2019estimating} for recent studies of the sample-splitting EM algorithm in mixture linear models. While the above-mentioned works all focus on low-dimensional models and unpenalized EM algorithm, regularized EM algorithm for high-dimensional mixture linear models is of growing interest in recent years. A selective review is as follows.
	
	\cite{khalili2007variable} studied the variable selection for mixture linear regression with penalized likelihood. When $p$ is fixed and $n$ goes to infinity, they established variable selection consistency and root-$n$ consistency for a possible local maximum. \citet{stadler2010L1} proposed a lasso-type penalty \citep{tibshirani1996regression} on the negative-log-likelihood function and obtained non-asymptotic convergence results for the global optimum. However, there is no guarantee for the EM algorithm to attain either the particular local maximum or the global maximum. To rigorously study the entire EM iterative solution sequence in high dimensions, existing results rely on the sample-splitting procedure. For example, convergence results based on sample-splitting algorithms are established for a truncated EM algorithm \citep{wang2014high}, a penalized EM algorithm \citep{yi2015regularized}, and a stochastic EM algorithm \citep{zhu2017high}, all under the assumptions of $K=2$, $\omega_1=\omega_2=1/2$, $\beta_2=-\beta_1$, and normally distributed predictors. More recently, with the help of sample splitting, \cite{zhang2020estimation} systematically studied estimation, confidence intervals, and large-scale hypotheses testing for the mixture linear model. Sample splitting is undoubtedly used to facilitate theoretical analyses, but is not desirable in practice.  It remains unknown how to practically choose $T$, which is both the number of iterations and the number of sample batches, and how it affects the estimation. Moreover, it is believed \citep[e.g.,][]{zhang2020estimation} that data splitting is unnecessary in numerical studies.

	The contributions of this article are multi-fold. The first and most significant contribution is
	developing a practical penalized EM algorithm for a general high-dimensional mixture linear model and establishing its substantial theoretical guarantees. Our algorithm is equally applicable to random-$X$ and fixed-$X$ and multiple mixtures $K\geq2$. In our theoretical analysis, we do not require sample splitting and allow a relatively general model. Specifically, we establish a non-asymptotic convergence rate for a two-mixture linear model with unknown proportions  $\omega_1,\omega_2\in(0,1)$, unrelated two regression parameter vectors  $\beta_1,\beta_2\in\mbbR^p$, and normally distributed predictors with unknown covariance structure. To our best knowledge, we established the first theoretical results for the high-dimensional mixture linear regression under such a general setting without sample splitting. Compared with the high-dimensional linear model literature, the theoretical analysis for the high-dimensional mixture model requires bounding the supremum of random processes and is much more challenging without sample splitting. Our general proof strategy is related to \cite{cai2019chime}, which studies the penalized EM algorithm for the Gaussian mixture model without data splitting. The complicated relationship between the random response $Y$ and the random predictor $X$ makes the theoretical studies of mixture linear regression even more challenging than the Gaussian mixture model, which only involves random $X$. Many new concentration results are needed for random processes involving both $Y$ and $X$. For instance, unlike $X_i$ in the Gaussian mixture model that is sub-Gaussian, the product term $X_iY_i$ appears frequently in the EM iterates and estimates and is more difficult to bound \citep{adamczak2008tail}. Even for the theoretical analysis of the population EM iterates, double expectations $\E_X\{\E_{Y\mid X}(\cdot)\}$ are needed than single expectation $\E_X(\cdot)$ in the Gaussian mixture model.
	With substantial efforts, we obtain a near optimal convergence rate of $\logg n\sqrt{s\logg p/n}$, with a small price $\logg(n)$ to pay for not sample splitting. 
	
	The second contribution is our new theoretical insights on model misspecification. Specifically, we analyze how a fixed parameter value $\sigma^2$ in the penalized EM algorithm may affect the estimation of $\beta$. To the best of our knowledge, it has not been studied in the literature on the mixture linear regression model. In most theoretical studies considering the EM algorithm for the mixture linear regression, the variance $\sigma^2$ is often assumed to be the true parameter value $(\sigma^{\star})^2$ and is a fixed constant. \cite{kwon2020minimax} considered the convergence results when $\sigma^2$ is updated in each M step. However, they only considered a simplified case where the mixture proportions are known to be $1/2$ and are not updated in each M step. Besides, their theoretical studies are limited to low dimensions. In general, when the error variance $\sigma^2$ is not correctly specified in the penalized EM algorithm brings bias for estimation. However, for the mixture model with a relatively large signal-to-noise ratio, our theory indicates that the choice of $\sigma^2$ has a minor influence on the estimation, and thus an accurate estimation for it is usually not necessary. This conclusion is further demonstrated by a simulation study.
	

	The third contribution is that we extend the study for the mixture linear regression model to multiple response cases. For the mixture linear regression model with a multivariate response, the naive approach is fitting a mixture linear regression model separately for each element of the response. The major drawback of doing so is each observation may be identified into different clusters when we model each univariate response separately. We illustrated the advantages of considering multiple responses together than handling them separately from both theoretical and numerical aspects. The advantage of considering multiple responses together is also demonstrated in \cite{hyun2020modeling}, which developed a
	sparse mixture linear regression model to estimate the time-varying data sets (i.e. at each time point, the data satisfy a mixture linear regression model). However, they are interested in developing a time-varying model to analyze a real-world dataset, but no theoretical study is conducted. In contrast, we rigorously characterize the advantage of considering multiple responses simultaneously with statistical theory.



The rest of the article is organized as follows. Section \ref{sec:est} contains the implementation details and discussions about the penalized EM algorithm. Section \ref{sec:the} presents the theory for the penalized EM algorithm and the influence of the choice of $\sigma^2$. Simulation studies are presented in Section \ref{sec:num}. We then extend the mixture linear regression to multiple response cases and consider its theoretical studies in Section \ref{MMLR}. In Section \ref{sec:real}, we consider a real data example followed by a short discussion in Section \ref{sec:dis}. The supplementary material contains proofs for all the Lemmas and Theorems and additional implementation details.

\section{Estimation}\label{sec:est}
We assume that we collect $n$ independent data points $\{(X_i,Y_i)\}_{i=1}^n$ from model \eqref{mlr}. In this section, we do not make additional assumptions on $X$ as our estimation procedure is equally applicable to random or fixed $X$ and allows for both continuous and discrete predictors. Let $\boltheta=\{ \omega_1,\cdots,\omega_k,\bolbeta_1,\cdots,\bolbeta_k\}$ be the unknown parameters to be estimated.
In this section, we focus on the studies of the regression coefficients and treat $\sigma^2$ as known. The estimation for $\sigma^2$ is briefly discussed at the end of this section. We further show in Theorem \ref{the:sigma} that misspecification of $\sigma^2$ has a relatively small impact on the final estimation.

To motivate our proposal, we first derive the standard EM algorithm and discuss its limitations.
The EM algorithm aims to maximize the log-likelihood of $Y\mid X$ over $\theta$, by iteratively updating the sequence of solutions $\{\widehat\theta^{(t)}, t=0,1,\ldots\}$ via the Expectation-step (E-step) and the Maximization-step (M-step). Recall that $W$ is the latent variable representing the mixtures. Consider the $(t+1)$-th iteration with the current value $\widehat{\theta}^{(t)}$. In the E-step, we calculate the expectation of the log-likelihood of $W\mid(Y,X)$ at the parameter $\widehat{\theta}^{(t)}$ and the observed $(X,Y)$. This is known as the Q-function,
\begin{equation}\label{Q}
\begin{aligned}
	Q(\theta\mid\widehat{\theta}^{(t)})&=-\frac{1}{2n}\sum_{i=1}^n\sum_{k=1}^K\widehat{\eta}_{i,k}(\widehat{\theta}^{(t)})(Y_i-X_i^T\beta_{k})^2 +\frac{1}{n}\sum_{i=1}^n\sum_{k=1}^K\widehat{\eta}_{i,k}(\widehat{\theta}^{(t)})\logg(\omega_k),
\end{aligned}
\end{equation} 
where $\widehat{\eta}_{i,k}(\widehat{\theta}^{(t)})=\mbbP(W_i=k\mid Y_i,X_i,\widehat{\theta}^{(t)} )$. The estimated probability $\widehat{\eta}_{i,k}(\widehat{\theta}^{(t)})$ is given by
\begin{equation}\label{eq:eta.raw}
\begin{aligned}
	\widehat{\eta}_{i,k}(\widehat{\theta}^{(t)})&=\frac{\widehat{ \omega}_k^{(t)}\phi_{\sigma^2}(Y_i-X_i^T\hatbolbeta_k^{(t)})}{\sum_{k=1}^K\widehat{ \omega}_k^{(t)}\phi_{\sigma^2}(Y_i-X_i^T\hatbolbeta_k^{(t)})},
\end{aligned}
\end{equation} 
where $\phi_{\sigma^2}(u)$ is the probability density function of $N(0,\sigma^2)$. Then, in the M-step, we update $\widehat{\theta}_k^{(t+1)}=\argmax_{\theta}Q(\theta\mid \widehat{\theta}^{(t)})$ by maximizing \eqref{Q}. 

Note that the standard EM algorithm is infeasible to high-dimensional problems. If $p\gg n$, even when the latent random variables $W_i$'s are observed, the maximizer of the Q function is not well-defined. Moreover, as $W_i$ is generally latent and unobserved, we need to calculate $\widehat{\eta}_{i,k}(\widehat{\theta}^{(t)})$ in \eqref{eq:eta.raw}, which involves the $p$-dimensional random vector $X$ and the $p$-dimensional parameter vectors $\beta_k$'s.

To estimate the linear mixture regression model in high dimensions, we modify the standard EM algorithm by encouraging sparsity. In high-dimensional statistics, it is often assumed that the coefficients have many elements as zero, i.e, most elements in $\bolbeta_k$ are zero. But we further assume that $\bolbeta_k$ has a joint sparsity structure, in that, for most $j$, we have $\beta_{1j}=\ldots=\beta_{Kj}=0$. The group sparsity facilitates the interpretation, as for each $j$, if $\beta_{1j}=\ldots=\beta_{Kj}=0$, then $X_j$ is unimportant for the prediction of $Y$ regardless of which mixture the observation comes from. Moreover, the group sparsity benefits the E-step, because, with straightforward calculation, we can rewrite \eqref{eq:eta.raw} as
\begin{equation}\label{eta}
\widehat{\eta}_{i,k}(\widehat{\theta}^{(t)})={\widehat{\omega}_k^{(t)}}/\big(\widehat{\omega}_k^{(t)}+\sum_{k'\neq k}\widehat{\omega}_{k'}^{(t)}\exp\big\{(\widehat{\beta}_{k'}^{(t)}-\widehat{\beta}_{k}^{(t)})^T X_i(Y_i-{(\widehat{\beta}_k^{(t)}+\widehat{\beta}_{k'}^{(t)})^TX_i}/{2})/\sigma^2  \big\} \big).
\end{equation}

Equation~\eqref{eta} shows an advantage of the group sparsity over the individual sparsity. It implies $X_j$ is unimportant for the evaluation of $\widehat{\eta}_{i,k}(\widehat{\theta}^{(t)})$ if $\widehat{\beta}_{k'}^{(t)}-\widehat{\beta}_{k}^{(t)}=0$ for all $k,k'$. The group sparsity guarantees that such a situation happens for most $X_j$ and $\widehat{\eta}_{i,k}(\widehat{\theta}^{(t)})$ is determined by a few elements in $\mbX$. 

With the sparsity assumption, we modify the EM algorithm by imposing the group lasso penalty \citep{yuan2006model} on $\beta_k$, for $k=1,\cdots,K$. Our the penalized EM algorithm replaces M-step by 
\begin{equation}\label{pQ}
\widehat{\theta}_k^{(t+1)}=\argmax_{\theta} \left\{2Q(\theta\mid \widehat{\theta}^{(t)}) - \lambda_n^{(t+1)} \sum_{j=1}^p\sqrt{\sum_{k=1}^K\beta_{kj}^2} \right\},
\end{equation}
where
$\lambda_n^{(t+1)}>0$ is the tuning parameter at the $(t+1)$-th iteration and $\beta_{kj}$ be the $j$-th element of $\beta_k$. In the E-step, we evaluate $\widehat{\theta}_k^{(t+1)}$ by \eqref{eta}, which is the same as in the standard EM algorithm. Clearly, our penalized EM algorithm reduces to the standard EM algorithm when $\lambda_n^{(t+1)}=0$ for all $t=0,1,\ldots$.

The optimization in \eqref{pQ} is separable in $\omega_k$ and $\bolbeta_k$. For $\omega$, the updating equation is the same as in the standard EM algorithm, i.e, $\widehat{\omega}_k^{(t+1)}=\sum_{i=1}^n\widehat{\eta}_{i,k}(\widehat{\theta}^{(t)})/n$. For $\beta_k$, it amounts to minimizing the following objective function, 
\begin{equation}\label{obj2}
\begin{aligned}
	\ell(\beta_1,\cdots,\beta_K)&=\sum_{k=1}^K\beta_k^T\widehat{\Sigma}_k^{(t+1)}\beta_k-2\sum_{k=1}^{K}(\widehat{\rho}_k^{(t+1)})^T\beta_k +\lambda_n^{(t+1)}\sum_{j=1}^p\sqrt{\sum_{k=1}^K\beta_{kj}^2},
\end{aligned}
\end{equation}
where $\widehat{\rho}_k^{(t+1)}=\sum_{i=1}^n\widehat{\eta}_{i,k}(\widehat{\theta}^{(t)})X_iY_i/n$ and $\widehat{\Sigma}_k^{(t+1)}= \sum_{i=1}^n\widehat{\eta}_{i,k}(\widehat{\theta}^{(t)})X_iX_i^T/n$. The convex optimization in \eqref{obj2} can be done efficiently by the groupwise majorization descent algorithm \citep{yang2015fast}. We provide implementation details in Supplement Section S.1.

The tuning parameter $\lambda_n^{(t)}$ could either be fixed or varying across iterations. For theoretical consideration, in Algorithm \ref{alg:mlr}, we set  $\lambda_n^{(t+1)}=\kappa\lambda_n^{(t)}+C_{\lambda}\sqrt{\logg(p)\logg(n)^2/n}$, where $0<\kappa<1/2$ and $C_{\lambda}$ are generic constants, for ease of showing the statistical convergence results. Note that $\lambda_n^{(t)}$ is at the order of $\sqrt{\logg(p)\logg(n)^2/n}$ when $t$ is large. Thus, in practice, we fix $\lambda_n^{(t)}=\lambda$ for all $t$ and tune $\lambda$ by the Bayesian information criterion (see our numerical studies). For fixed $\lambda_n^{(t)}=\lambda$ over all iterations, our penalized EM algorithm is maximizing
\begin{equation}\label{pl}
L(\theta)-\lambda/2 \sum_{j=1}^p\sqrt{\sum_{k=1}^K\beta_{kj}^2},
\end{equation}
where $L(\theta)$ is the conditional log-likelihood of $Y\mid X$. The following lemma shows the convergence result of Algorithm \ref{alg:mlr}.

\begin{lemma}
If we set $\lambda_n^{(t)}=\lambda$ for all $t$ in Algorithm \ref{alg:mlr}, the objective function from \eqref{pl} evaluated at $\theta^{(t+1)}$ is guaranteed to be no less than
the objective function from \eqref{pl} evaluated at $\theta^{(t)}$. That is, the sequence of iterates $\{\theta^{(t)}\}_{t=1}^{\infty}$ generated by Algorithm \ref{alg:mlr} monotonically increase the value of the objective function from \eqref{pl}.
\end{lemma}

In Algorithm \ref{alg:mlr}, $\sigma^2$ is treated as a known parameter to facilitate theoretical studies. Treating $\sigma^2$ as known is also common in theoretical studies for mixture linear regression \citep[e.g.]{yi2015regularized, balakrishnan2017statistical, zhang2020estimation}. We leave the detailed discussion for it in Section \ref{sec:the}.
In practice, the estimates of $\sigma^2$ can be updated straightforwardly in the EM algorithm as $n^{-1}\sum_{i=1}^n\sum_{k=1}^K\widehat{\eta}_{i,k}(\widehat{\theta}^{(t+1)})(Y_i-X_i^T\widehat{\beta}_k^{(t+1)})^2$. Similar estimates are also adopted in numerical studies of \cite{zhang2020estimation}.

Regularization strategies are also used by \cite{yi2015regularized}, \cite{cai2019chime}, and \cite{zhang2020estimation} in high-dimensional EM algorithms. However, \cite{cai2019chime} considers the classification problem instead of regression problem. \cite{yi2015regularized,zhang2020estimation} studies the regression problem with the addition of the lasso penalty \citep{tibshirani1996regression} instead of the group lasso penalty. However, a key difference between our algorithm and theirs is that they require data to be split into $T$ batches, and the (penalized) EM algorithm iterates $T$ times, using one independent batch at each iteration. The sample splitting is rarely performed in standard EM algorithms on low-dimensional data, as it may decrease the computation efficiency. Instead, the sample splitting is an attempt to circumvent technical difficulty in proving the convergence rate. Our proposed EM algorithm does not split sample, but we will show that it achieves a high level of accuracy regardless. Moreover, in addition to investigating the property of the penalized EM algorithm in estimating model \eqref{mlr}, we also study the effect of misspecification of $\sigma^2$ and the estimation of mixture linear regression when there are multiple responses; see Theorem~\ref{the:sigma} and Section~\ref{MMLR}, respectively.


\begin{algorithm}[t!]
\caption{Group lasso penalized EM algorithm for model \eqref{mlr} in high dimensions}
\label{alg:mlr}
\begin{enumerate}
	\item[Input:] Initial values $\widehat{\omega}_k^{(0)}$, $\widehat{\beta}_k^{(0)}$, for $k=1,\cdots,K$, maximum iteration number $T$, data $\{X_i,Y_i; i=1,\ldots,n\}$, and initial tuning parameter	
	\begin{equation*}
		\begin{aligned}
			\lambda_n^{(0)}=C_1(\vert\widehat{ \omega}_1^{(0)}-\omega_1^*\vert\lor\Vert\widehat{\beta}_1^{(0)}-\beta_1^*\Vert_2\lor\cdots\lor\quad\Vert\widehat{\beta}_K^{(0)}-\beta_K^*\Vert_2)/{\sqrt{s}}+C_{\lambda}\sqrt{\logg(n)^2\logg(p)/s}
		\end{aligned}
	\end{equation*}
	for some positive constants $C_1$ and $C_{\lambda}$. 
	
	\item[Iterate:] {For $t=0,\cdots,T-1$, do the following steps until convergence.} 
	\begin{itemize}
		
		\item For $i=1,\cdots,n$, let
		\begin{equation*}
			\begin{aligned}
				&\quad \widehat{\eta}_{i,k}(\widehat{\theta}^{(t)})={\widehat{\omega}_k^{(t)}}/\Big(\widehat{\omega}_k^{(t)}+\sum_{k'\neq k}\widehat{\omega}_{k'}^{(t)}\exp\big\{(\widehat{\beta}_{k'}^{(t)}-\widehat{\beta}_{k}^{(t)})^T\cdot X_i\big(Y_i-{(\widehat{\beta}_k^{(t)}+\widehat{\beta}_{k'}^{(t)})^TX_i}/{2}\big)/\sigma^2  \big\} \Big).
			\end{aligned}
		\end{equation*}
		
		\item For $k=1,\cdots, K$, update 
		\begin{equation*}
			\begin{aligned}
				&\widehat{\omega}_k^{(t+1)}=\frac{1}{n}\sum_{i=1}^n\widehat{\eta}_{i,k}(\widehat{\theta}^{(t)}),\\ &\widehat{\rho}_k^{(t+1)}=\frac{1}{n}\sum_{i=1}^n\big(\widehat{\eta}_{i,k}(\widehat{\theta}^{(t)})X_iY_i\big),\\ &\widehat{\Sigma}_k^{(t+1)}=\frac{1}{n}\sum_{i=1}^n\big(\widehat{\eta}_{i,k}(\widehat{\theta}^{(t)})X_iX_i^T\big),
			\end{aligned}
		\end{equation*}
		and update $\widehat{\beta}_k^{(t+1)}$ by minimizing
		\begin{equation*}
			\begin{aligned}
				\sum_{k=1}^K\beta_k^T\widehat{\Sigma}_k^{(t+1)}\beta_k-2\sum_{k=1}^{K}(\widehat{\rho}_k^{(t+1)})^T\beta_k+\lambda_n^{(t+1)}\sum_{j=1}^p\sqrt{\sum_{k=1}^K\beta_{kj}^2}.
			\end{aligned}
		\end{equation*}
		with $\lambda_n^{(t+1)}=\kappa\lambda_n^{(t)}+C_{\lambda}\sqrt{\logg(p)\logg(n)^2/n}$, where $\kappa\in(0,1/2)$.
	\end{itemize}
	\item[Output:] $\widehat\theta^{(t+1)}=\{\widehat\omega_1^{(t+1)},\ldots,\widehat\omega_K^{(t+1)},\widehat\beta_1^{(t+1)},\ldots,\widehat\beta_K^{(t+1)}\}$
\end{enumerate}
\end{algorithm}

\section{Theory}\label{sec:the}
\subsection{Preliminary}
We begin this section with some notations. For numbers $a$ and $b$, $a\lor b$ means $\max\{a, b\}$. For an integer $n$, we let $[n]$ denote the set $\{1,\cdots,n \}$. For a vector $x=(x_1,\cdots,x_p)^T$, $\Vert x\Vert_0$ is the number of non-zero elements in $x$, $\Vert x\Vert_1=\sum_{i=1}^{p}\vert x_i\vert$, and $\Vert x\Vert_2=\sqrt{\sum_{i=1}^p x_i^2}$.
For a symmetric matrix $A$, we denote $\lambda_{min}(A)$ and $\lambda_{max}(A)$ as the smallest and largest eigenvalues of $A$, respectively. The Frobenius norm of a matrix $A=(a_{ij})$ is  defined as $\Vert A\Vert_F=\sqrt{\sum_{i,j}a_{ij}^2}$. The $\ell_2$ norm of a matrix $A$ is  $\Vert A\Vert_2=\sqrt{\lambda_{max}(A^TA)}$.  For a subset $\calA\subseteq \{1,\cdots,p\}$, $\calA^c$ denotes its complement. For two sequences of positive numbers $a_n$ and $b_n$, $a_n=O(b_n)$ means $a_n\leq cb_n$ for a constant $c>0$ for all $n$, $a_n=o(b_n)$ means that $a_n/b_n\to 0$ as $n\to \infty$, and $b_n\gg a_n$ means that $a_n=o(b_n)$. Let $\calS^{p-1}$ be the unit sphere. For a positive integer $s\leq p/2$, let set $\Gamma_{2p}(2s)=\{\mu\in\mbbR^{2p}: \Vert\mu_{S^c}\Vert_1\leq 5\sqrt{2s}\Vert\mu_{S}\Vert_2+2\sqrt{2s}\Vert\mu\Vert_2$ for some $S\subset [2p]$ with $\vert S\vert=4s\}$ and $\Gamma(s)=\Gamma_{2p}(2s)_{1:p}$, where $\Gamma_{2p}(2s)_{1:p}=\{\mu_{1:p}:\mu\in\Gamma_{2p}(2s)  \}$. 
For a vector $x$ and a symmetric matrix $A$, we define $\Vert x\Vert_{2,s}=\sup_{\Vert\mu\Vert_2=1,\mu\in\Gamma(s)}\langle x,\mu\rangle$, and $\Vert A\Vert_{2,s}=\sup_{\Vert\mu\Vert_2=1,\mu\in\Gamma(s)}\vert\mu^TA\mu\vert$. 

We assume independent random predictors $X_i\sim N(0,\Sigma)$ in our theoretical analysis. This is less restrictive than the assumption of $X_i\sim N(0,I_p)$ in \cite{yi2015regularized} and \cite{balakrishnan2017statistical} in that we allow the predictor to be correlated. In this section, we consider $K=2$, which is a common assumption in theoretical analysis for high-dimensional EM algorithm \citep[e.g.]{yi2015regularized,cai2019chime,zhang2020estimation}.  We re-define $\theta=\{\omega_1,\beta_1,\beta_2\}$ since $\omega_2=1-\omega_1$. Let $\boltheta^*$ be the true value of $\boltheta$, and $\hatboltheta^{(t)}$ be the estimate of $\boltheta$ at the $t$-th step in Algorithm~\ref{alg:mlr}. 
The true parameter space we consider is 
\begin{equation*}
\begin{aligned}
	\bolTheta^*=\{\boltheta^*: \omega_1^*\in (c_w,1-c_w), \Vert\bolbeta_k^*\Vert_0\leq s, \Vert\bolbeta_k^*\Vert_2\leq M_b,\ \mathrm{for}\ k=1,2 \}. 
\end{aligned}
\end{equation*}
This is a natural parameter space to consider. The condition $\omega_1^*\in (c_w,1-c_w)$ guarantees the sample size from each latent class is large enough. Condition on $\Vert\bolbeta_k^*\Vert_2\leq M_b$ is also similarly used in \cite{yi2015regularized} under the data splitting framework and is milder than a similar condition $\Vert\bolbeta\Vert_1\leq M_b$ used in \cite{cai2019chime}, where the EM algorithm for Gaussian mixture model without data splitting is studied. 

In theoretical studies, we first assume that the true value of $\sigma^2$ denoted as $(\sigma^{\star})^2$ as a known parameter, namely, the input $\sigma^2$ is $(\sigma^{\star})^2$ in Algorithm \ref{mlr}.
Without loss of generality, we assume $(\sigma^{\star})^2=1$. Treating $\sigma^2$ as known is also common in state-of-the-art theoretical studies for mixture linear regression \citep[e.g.]{yi2015regularized, balakrishnan2017statistical, zhang2020estimation}. Although we do not analyze $\bolbeta_k$'s and $\sigma^2$ simultaneously,  we investigate how the choice of $\sigma^2$ influence the estimation of Algorithm \ref{alg:mlr} in Theorem \ref{the:sigma}.

Since $\sigma^2=1$ and $K=2$, we simplify $\mbbP(W_i=1\mid Y_i,X_i,{\theta} )$ as
\begin{equation*}
\begin{aligned}
	\quad \eta_{i,1}(\theta)=1-\eta_{i,2}(\theta)=1/\big[1+(\omega_2/\omega_1)\exp\{(\bolbeta_2-\bolbeta_1)^T\mbX_i\cdot(Y_i-{(\bolbeta_1+\bolbeta_2)^T\mbX_i}/{2}) \}\big].
\end{aligned}
\end{equation*}
Then the following quantities are used repeatedly in our theoretical analysis:
\begin{equation*}
\begin{aligned}
	&\widehat{\omega}_k(\boltheta)=\frac{1}{n}\sum_{i=1}^n\eta_{i,k}(\theta),\
	\widehat{\bolrho}_k(\boltheta)=\frac{1}{n}\sum_{i=1}^n\eta_{i,k}(\theta)\mbX_iY_i,\\
	&\hatbolSigma_k(\boltheta)=\frac{1}{n}\sum_{i=1}^n\eta_{i,k}(\theta)\mbX_i\mbX_i^T,\
	{\omega}_k(\boltheta)=\E\big\{\frac{1}{n}\sum_{i=1}^n\eta_{i,k}(\theta)\big\},\\
	&{\bolrho_k}(\boltheta)=\E\big\{\frac{1}{n}\sum_{i=1}^n\eta_{i,k}(\theta)\mbX_iY_i\big\},\
	\bolSigma_k(\boltheta)=\E\big\{\frac{1}{n}\sum_{i=1}^n\eta_{i,k}(\theta)\mbX_i\mbX_i^T\big\},
\end{aligned}
\end{equation*}
where the expectation is with respect to $\mbX_i$ and $Y_i$, $i=1,\ldots,n$.  
Then we let $M(\boltheta)=\{{\omega}_k(\boltheta),{\bolrho_k}(\boltheta),\bolSigma_k(\boltheta), k=1,2 \}$, $M_n(\boltheta)=\{\widehat{\omega}_k(\boltheta),\widehat{\bolrho}_k(\boltheta),\hatbolSigma_k(\boltheta), k=1,2\}$, and define $d_{2,s}\big( M(\boltheta_1),M(\boltheta_2)\big)$ and $d_{2}\big( M(\boltheta_1),M(\boltheta_2)\big)$ as
\begin{equation*}
\begin{aligned}
	&\max_{k=1,2}\{\vert\omega_k(\boltheta_1)-\omega_k(\boltheta_2)\vert \lor\Vert\bolrho_k(\boltheta_1)-\bolrho_k(\boltheta_2)\Vert_{2,s}\lor\Vert(\bolSigma_k(\boltheta_1)-\bolSigma_k(\boltheta_2))\bolbeta_2^*\Vert_{2,s}\},\\
	&\max_{k=1,2}\{\vert\omega_k(\boltheta_1)-\omega_k(\boltheta_2)\vert \lor\Vert\bolrho_k(\boltheta_1)-\bolrho_k(\boltheta_2)\Vert_{2}\lor\Vert(\bolSigma_k(\boltheta_1)-\bolSigma_k(\boltheta_2))\bolbeta_2^*\Vert_{2}\},
\end{aligned}
\end{equation*}
respectively, which are distances between $M(\boltheta_1)$ and $M(\boltheta_2)$.

Let $\Delta=\sqrt{(\bolbeta_2^*-\bolbeta_1^*)^T\bolSigma (\bolbeta_2^*-\bolbeta_1^*)}$, which is a measure of the signal-to-noise ratio of the mixture  linear regression model.
We define the contraction basin $\calB_{con}(\boltheta)$ as follows.
\begin{equation*}
\begin{aligned}
	\calB_{con}(\boltheta)=\{\boltheta:\omega_k\in (c_0,1-c_0), \Vert\bolbeta_k-\bolbeta_k^*\Vert_2\leq C_b\Delta,\  \bolbeta_k-\bolbeta_k^*\in\Gamma(s),  \ \mathrm{for}\ k=1,2\}.
\end{aligned}
\end{equation*}

Intuitively, the contraction basin requires that $\bolbeta_k$ is not far away from the true parameter $\bolbeta_k^*$. Under the technical conditions shown later, an initialization $\hatboltheta^{(0)}$ falls in the contraction basin can guarantee the subsequent estimators  $\hatboltheta^{(t)}$ in Algorithm \ref{alg:mlr} are all contained in the contraction basin. 


\subsection{Main results}
We first introduce some technical conditions before stating the theoretical results.
\begin{enumerate}[(C1)]
\item The eigenvalues of $\bolSigma$ satisfy that $M_1\leq \lambda_{min}(\bolSigma)\leq \lambda_{max}(\bolSigma)\leq M_2$.
\item $n\gg s\log(p)$.
\item	The signal-to-noise ratio $\Delta>C_1(c_0)$ for a constant $C_1(c_0)$ only depends on $c_0$, and $C_b<C_2(c_0,M_2)$ for a constant $C_2(c_0,M_2)$ only depends on $c_0,M_2$.
\item  The initialization $\hatboltheta_0^{(0)}=(\widehat{ \omega}_1^{(0)},\hatbolbeta_1^{(0)},\hatbolbeta_2^{(0)})\in\calB_{con}(\boltheta^*)$.
\end{enumerate}
Condition (C1) is a standard assumption on the covariance $\bolSigma$ in high-dimensional statistics \citep{bickel2008covariance, cai2011constrained}. Condition (C2) is a common assumption in high dimensions  on the relationship among $(n,p,s)$ to guarantee consistent estimation \citep[e.g]{meinshausen2009lasso}. In particular, it implies that the restrictive eigenvalue condition  $\inf_{\mu\in\Gamma(s)\cap\calS^{p-1}}\{\mu^T(\sum_{i=1}X_iX_i^T/n)\mu\}>\tau_0$ holds for a positive generic constant $\tau_0$ with high probability, and is used for proving the concentration of $\hatbolbeta_k^{(t)}$ in the $t$-th iteration.  Condition (C3) has two requirements. The first one is that the signal-to-noise ratio is larger than a universal constant that does not depend on $n$ and $p$ so that the two mixtures are distinguishable. This requirement was also previously used in mixture linear model \citep[e.g.,][]{yi2015regularized, balakrishnan2017statistical, zhang2020estimation}. The second one is that,  for the parameter $\bolbeta_k$ in the contraction basin, the distance $\Vert\bolbeta_k-\bolbeta_k^*\Vert_2$ is bounded by the signal-to-noise ratio multiplied by a generic constant independent of $n$ and $p$. This requirement makes all the $\bolbeta_k$ in the contraction basin not too far away from the truth $\bolbeta_k^*$.  Condition (C4) ensures that the initialization is in the contraction basin. The contraction and concentration properties shown later guarantee that the estimates in each step of Algorithm \ref{alg:mlr} stays in the contraction basin.  

Next, we present two lemmas about the linear convergence of the
population EM updates and the concentration of the sample estimation to the population one in each EM iteration. The following two lemmas together with a key technical Lemma S.10 in the Supplement are highly non-trivial and serve as the building blocks of the main theory for Algorithm \ref{alg:mlr}. 
\begin{lemma}\label{contraction}
Under conditions (C1) and (C3), if $\theta\in\calB_{con}(\boltheta^*)$, then
\begin{equation*}
	d_{2}\big(M(\boltheta),M(\boltheta^*)\big)\leq \kappa_0\big( \vert\omega_1(\boltheta)-\omega_1^*\vert\lor\Vert\bolbeta_1-\bolbeta_1^*\Vert_2\lor\Vert\bolbeta_2-\bolbeta_2^*\Vert_2 \big).
\end{equation*}
for some $0<\kappa_0<\frac{1}{2\lor (64/\tau_0)}$.
\end{lemma}
\begin{lemma}\label{concentation}
Suppose that $\boltheta^*\in\bolTheta^*$. Under condition (C1), there exists a constant $C_{con}>0$, such that with probability at least $1-4p^{-1}$,
\begin{equation*}
	\begin{aligned}
		&\sup_{\boltheta\in \calB_{con}(\boltheta^*)}d_{2,s}(M(\boltheta),M_n(\boltheta))\leq C_{con} \sqrt{\frac{s\logg(n)^2\logg(p)}{n}}
	\end{aligned}
\end{equation*}
\end{lemma}

Intuitively, Lemma \ref{contraction} shows the computational contraction of Algorithm \ref{alg:mlr}. It implies that, in each EM iterations, the expectations of the updated estimators converge to the true parameters at a linear rate. On the other hand, Lemma \ref{concentation} establishes the statistical convergence rate of estimators to their expectations in each EM iteration. When the iteration steps are large enough, the computational error will be dominated by the statistical error, which means that further iterations can not improve the statistical convergence rate of the algorithm. \cite{cai2019chime} also proved similar lemmas under the Gaussian mixture model, but our proof is more challenging, as we are interested in the mixture linear regression model. The unboundedness in both $X$ and $Y$ makes $M(\theta)$ more complicated and $M_n(\theta)$ have heavier tails. 

Thanks to Lemmas \ref{contraction} and \ref{concentation}, we can show the following result for Algorithm \ref{alg:mlr}.
\begin{theorem}\label{the:1}
Under conditions (C1)--(C4), there exists a constant $0<\kappa<1/2$, such that $\hatbolbeta_k^{(t+1)}$ obtained by Algorithm \ref{alg:mlr} satisfies, with probability $1-4p^{-1}$,
\begin{equation*}
	\begin{aligned}
		\Vert\hatbolbeta_k^{(t+1)}-\bolbeta_k^*\Vert_2&=O\Big(\kappa^t(\vert\widehat{ \omega}_1^{(0)}-\omega_1^*\vert\lor\Vert\hatbolbeta_1^{(0)}-\bolbeta_1^*\Vert_2\lor\Vert\hatbolbeta_2^{(0)}-\bolbeta_2^*\Vert_2) +\sqrt{\frac{s\logg(n)^2\logg(p)}{n}}\Big).
	\end{aligned}
\end{equation*}
Consequently, for $t\geq \{-\logg(\kappa)\}^{-1}\logg\{n(\vert\widehat{ \omega}_1^{(0)}-\omega_1^*\vert\lor\Vert\hatbolbeta_1^{(0)}-\bolbeta_1^*\Vert_2\lor\Vert\hatbolbeta_2^{(0)}-\bolbeta_2^*\Vert_2)\}$,
\begin{equation*}
	\Vert\hatbolbeta_k^{(t+1)}-\bolbeta_k^*\Vert_2=O\Big(\sqrt{\frac{s\log(n)^2\logg(p)}{n}}\Big).
\end{equation*}
\end{theorem}
We make several remarks on Theorem \ref{the:1}. Firstly,  compared with existing results \citep{yi2014alternating, yi2015regularized, balakrishnan2017statistical} requiring $\bolbeta_2=-\bolbeta_1$ and $\mbX\sim N(0,I_p)$, our model settings are more general. On one hand, note that $\bolbeta_2=-\bolbeta_1$ is not just a location shift of the response
variable. As an illustration, when $\bolbeta_1=(1,1,0)^T$ and $\bolbeta_2=(0.5,2,0)^T$, a simple location shift cannot reduce the model to the case where $\bolbeta_2=-\bolbeta_1$. By removing the assumption that $\bolbeta_2=-\bolbeta_1$, our theory is applicable to a larger model space. On the other hand, the predictors are not likely to be uncorrelated in practice. Thus, it is meaningful to extend the condition $\mbX\sim N(0,I_p)$ to $\mbX\sim N(0,\bolSigma)$. Also, although we adopt the group lasso penalty in the algorithm, the theoretical results can be naturally extended to the penalties that are decomposable \citep{negahban2012unified}, such as the lasso penalty popular in the literature.

Moreover, unlike the extensive literature about the sample-splitting EM algorithms for mixture  linear regression \citep{yi2014alternating,yi2015regularized,zhang2020estimation}, the convergence result in Theorem~\ref{the:1} does not require sample splitting. Sample splitting is not desirable, especially when we have a small sample size. Splitting a limited number of observations into $T$ batches decreases the estimation efficiency, makes the estimation less stable and is rarely used in practice. Hence, it is meaningful to develop theoretical results without data splitting for mixture linear regression, as those in Theorem~\ref{the:1}. To our best knowledge, Theorem \ref{the:1} is the first theoretical result for the high-dimensional EM algorithm of mixture linear regression without data splitting. Also note that  the convergence rate we obtain is nearly optimal. When the latent random variables $W_i$'s are unknown, the optimal rate $\sqrt{s\logg(p)/n}$ \citep[e.g.]{ye2010rate}, while when $W_i$'s are unknown, \citet{zhang2020estimation} gives an estimation rate of $\sqrt{s\logg(p)\logg{n}/n}$ with sample splitting. Our result is slightly slower than these rates by the factors of $\logg{n}$ and $\logg^{1/2}{n}$, respectively. The additional $\logg(n)$ terms are the price of no data splitting. Technically, to prove the convergence results without data splitting we have to bound the tails of the supremum of unbounded random processes, which is much more challenging than bounding the random variables.  

Now we turn to the effect of the misspecification of $\sigma^2$. In most existing theoretical analysis including this one, $\sigma^2$ is usually treated as a known parameter in theoretical studies for the mixture linear regression model. There are two reasons for this treatment. One one hand, the regression coefficients are of primary interest in regression models. On the other, statistical analysis for mixture regression model with known $\sigma^2$  is already challenging. However, in practice $\sigma^2$ is almost never known. Statisticians often plug in an estimated value of $\sigma^2$, which differs from $\sigma^2$. Yet it is unclear how the misspecification affects the estimation of the mixture linear regression. In the following theorem, we obtain a non-asymptotic convergence result for Algorithm \ref{alg:mlr} with misspecified $\sigma^2$, which indicates that although misspecified  $\sigma^2$ brings bias to the estimation, it usually has a minor influence on the mixture linear regression model with a large signal-to-noise ratio.
\begin{theorem}\label{the:sigma}
Let $\widetilde{\bolbeta}_k^{(t+1)}$ be the estimation of Algorithm \ref{alg:mlr} in the $(t+1)$-th EM step where $\sigma^2$ may not equal to $(\sigma^{\star})^2$. Under Conditions (C1)-(C4) and $ (\sigma^{\star})^2/2<\sigma^2<2(\sigma^{\star})^2$, with probability $1-4p^{-1}$,
\begin{equation*}
	\begin{aligned}
		\Vert\widetilde{\bolbeta}_k^{(t+1)}-\bolbeta_k^*\Vert_2&=O\Big(\kappa^t(\vert\widehat{ \omega}_1^{(0)}-\omega_1^*\vert\lor\Vert\hatbolbeta_1^{(0)}-\bolbeta_1^*\Vert_2\lor\Vert\hatbolbeta_2^{(0)}-\bolbeta_2^*\Vert_2)+\sigma^{\star}(\Delta/\sigma^{\star})^{-2}\vert 1-\frac{(\sigma^{\star})^2}{\sigma^2}\vert\\& +\sigma^{\star}\sqrt{\frac{s\logg(n)^2\logg(p)}{n}}\Big).
	\end{aligned}
\end{equation*}
Consequently, for $t\geq \{-\logg(\kappa)\}^{-1}\logg\{n(\vert\widehat{ \omega}_1^{(0)}-\omega_1^*\vert\lor\Vert\hatbolbeta_1^{(0)}-\bolbeta_1^*\Vert_2\lor\Vert\hatbolbeta_2^{(0)}-\bolbeta_2^*\Vert_2)\}$,
\begin{equation*}
	\Vert \widetilde{\bolbeta}_k^{(t+1)}-\beta_k^*\Vert_2\leq c_1\sigma^{\star}(\Delta/\sigma^{\star})^{-2}\vert 1-\frac{(\sigma^{\star})^2}{\sigma^2}\vert+c_2\sigma^{\star}\sqrt{\frac{s\logg (n)^2\logg p}{n}}.
\end{equation*}
\end{theorem}

Theorem \ref{the:sigma} considered the case where $(\sigma^{\star})^2$ is misspecified as $\sigma^2$ in Algorithm \ref{alg:mlr}.
The only difference between the convergence rates in Theorems \ref{the:1} and \ref{the:sigma} is the bias term $\sigma^{\star}(\Delta/\sigma^{\star})^{-2}\vert 1-{(\sigma^{\star})^2}/{\sigma^2}\vert$. It is a linear function of $|1-(\sigma^{\star})^2/\sigma^2|$ and is exactly 0 when $\sigma^2=(\sigma^{\star})^2$. Although not using $(\sigma^{\star})^2$ in Algorithm \ref{alg:mlr} may return a biased estimation, the bias is small when $\sigma^2$ is close to $(\sigma^{\star})^2$. Note that the bias term is also a function of $(\Delta/\sigma^{\star})^{-2}$. When the signal-to-noise ratio $\Delta$ is large, the bias term is small regardless of the value of $\sigma^2$. Specifically, if $\Delta\gg O( n^{1/4}/[s\logg(n)^2\logg p]^{1/4})$, then the bias term is ignorable compared with $\sqrt{s\logg(n)^2\logg p/n}$. For large $p$, small $n$ cases, the requirement $\Delta\gg O( n^{1/4}/[s\logg(n)^2\logg p]^{1/4})$ is usually not strict. Intuitively, when $\Delta\to \infty$, the mixtures can be easily identified, which makes the mixture linear regression model reduce to linear regression models. It is well known that the choice of $\sigma^2$ in the linear regression model does not affect the maximum likelihood estimation for the coefficients. Thus, when $\Delta\to\infty$, the bias term disappears.  Our numerical studies also demonstrate that $\sigma^2$ has a minor influence on the performance of Algorithm \ref{alg:mlr} for relatively large $\Delta$. Please see Section \ref{addsimu} for details. We remark that the condition $(\sigma^{\star})^2/2<\sigma^2<2(\sigma^{\star})^2$ is not necessary for the algorithm in practice, but only a requirement in theory due to technical reasons.

\section{Simulation studies}\label{sec:num}

\subsection{Simulation set-up}
In this section, we investigate the empirical performance of Algorithm~\ref{alg:mlr}. For practical initialization, we start with lasso regression \citep{tibshirani1996regression} of $Y$ on $X$ to roughly select important predictors and then apply the tensor power method \citep{anandkumar2014tensor} on the selected variables to get initializations for $\bolbeta_k$ and mixing proportion $\omega_k$, $k=1,\cdots,K$. We use $T=20$ as the maximum number of iterations in Algorithm~\ref{alg:mlr}, and stop the iterations when $\Vert(\widehat\beta^{(t+1)}_1,\ldots,\widehat\beta^{(t+1)}_K)-(\widehat\beta^{(t)}_1,\ldots,\widehat\beta^{(t)}_K)\Vert_F<10^{-3}$ . We use a single tuning parameter $\lambda=\lambda_n^{(t)}$ for all $t=0,1,\dots,T$ and choose $\lambda$ based on the Bayesian information criterion. 

We include the following methods: 1) \emph{Oracle}, we fit the standard EM algorithm on the true subset of $s$ relevant predictors. 2) \emph{Initial}, after using the tensor power method to get the initialization for $\beta_k$, $k=1,\cdots,K$, we plug them into \eqref{eta} to get an estimation for the weight $\widehat{\eta}_{ik}$ to the $i$-th sample, $i=1,\cdots,n$. Then, we label the $i$-th sample with $\argmax_{k=1,\cdots,K}\widehat{\eta}_{ik}$. After that, we fit lasso regressions separately for each estimated mixtures to estimate $\beta_k$. Compared with the tensor power method, this method returns sparse estimations. 3) \emph{GLLiM}, the Gaussian Local Linear Mapping EM algorithm \citep{deleforge2015high} that is implemented in the R package \texttt{xLLiM}. 4) \emph{HDEM}, our implementation of the high-dimensional EM algorithm proposed by \cite{zhang2020estimation}. 5) \emph{PSEM}, i.e., the post-selection EM algorithm, we fit the standard EM algorithm on the selected variables from Algorithm~\ref{alg:mlr}. and finally 6) \emph{PEM}, the penalized EM algorithm (Algorithm~\ref{alg:mlr}).

We consider the following simulation models, where we first generate independent  predictor $X_i\sim N(0,\Sigma)$ and error $\epsilon_i\sim N(0,1)$ and then $Y_i$ follows from the mixture linear regression \eqref{mlr}. For all the four simulation examples, we fix the first $s=10$ coefficients in $\beta_k$ to be nonzero and consider both $p=400$ and $p=1000$ settings. The total sample size $n$ is set to be $400$ for models (M1)--(M3), where we have two mixtures, and $600$ for model (M4), where we have three mixtures. The symmetric mixture assumption does not hold in any of our models, $\beta_1\neq -\beta_2$.
\begin{itemize}
\item (M1) Two mixtures with $\omega_1=\omega_2=0.5$, and auto-regressive covariance structure $[\bolSigma]_{ij}=0.3^{\vert i-j\vert}$ for $i, j=1,\ldots,p$. The nonzero coefficients of $\bolbeta_1$ is generated independently from $N(0,1)$; and the nonzero coefficients of $\bolbeta_2$ is $\bolbeta_{2j}=\bolbeta_{1j}+2\cdot\mathrm{sgn}(\bolbeta_{1j})$, $j=1,\ldots,s$.
\item (M2) Same as (M1) but with weaker signals:$\bolbeta_{2j}=\bolbeta_{1j}+1\cdot\mathrm{sgn}(\bolbeta_{1j})$, $j=1,\ldots,s$.
\item (M3) Same as (M1) but with a different covariance $\Sigma$ based on Erd\'os-R\'enyi random graph. Let $\tilbolSigma=(\widetilde{\sigma}_{ij})$, where $\widetilde{\sigma}_{ij}=u_{ij}\delta_{ij}$, $\delta_{ij}$ follows Bernoulli$(0.1)$ distribution, and $u_{ij}\sim \mathrm{Uniform}[0.5,1]\cup  \mathrm{Uniform}[-1,-0.5]$. Then let $\tilbolSigma_1=(\tilbolSigma+\tilbolSigma^T)/2$ and $\bolSigma^*= \tilbolSigma_1+\{\max(-\lambda_{min}(\tilbolSigma_1),0 )+0.05   \}I_p$. Finally, $\bolSigma^*$ is standardized to have $1$'s on the diagonal.
\item (M4) Three mixtures with $\omega_1=\omega_2=\omega_3=1/3$. The nonzero elements of $\bolbeta_1$ and $\beta_3$ are $-1$ and $5$, respectively, and the nonzero elements of $\bolbeta_2$ are evenly spaced between $1$ and $3$.
\end{itemize}

\subsection{Simulation results}

The performances of those methods are evaluated using the following criteria. The parameter estimation errors are defined as $\sqrt{\sum_{k=1}^K\Vert\hatbolbeta_k-\bolbeta_k\Vert_2^2}$ for $\beta$ and $\sum_{k=1}^K\vert\widehat{\omega}_k-\omega_k\vert\times 100$ for $\omega$. The mixture estimation error is defined as $\sum_{i=1}^n I(\widehat{W}_i\neq W_i)/n\times 100$, where $\widehat W_i = \argmax_{k=1,\ldots,K} \widehat\eta_{i,k}(\widehat\theta)$ for $i=1,\dots,n$. For methods that performs variable selection, we also recorded the true positive and false positive rates of variable selection.
%
%
\begin{table*}[ht!]
\centering
\resizebox{\textwidth}{!}{
	\renewcommand\arraystretch{0.91}
	\begin{tabular}{ccccccc}
		\hline 			
		\multirow{2}{*}{\textbf{M1} } & \multicolumn{3}{c}{$p=400$} & \multicolumn{3}{c}{$p=1000$}\tabularnewline
		& $\beta_k$ & $\omega_k$ & $W_i$ 	& $\beta_k$ & $\omega_k$ & $W_i$  \tabularnewline
		\hline 
		\multirow{1}{*}{Oracle }& 0.44 (0.004) & 6.08 (0.06) & 8.43 (0.01)  & 0.39 (0.001) & 6.18 (0.05) & 8.40 (0.01) \tabularnewline 
		\multirow{1}{*}{PSEM }& 0.57 (0.01) & 6.31 (0.06) & 9.17 (0.04)  & 0.76 (0.01) & 6.28 (0.05) & 10.50 (0.08) \tabularnewline               
		\multirow{1}{*}{GLLiM }& 9.63 (0.01) &68.18 (0.23)  & 33.23 (0.12)  & 9.68 (0.01) & 74.20 (0.23) & 34.24 (0.13) \tabularnewline
		\multirow{1}{*}{Initial }& 5.14 (0.01) & 32.14 (0.20) & 41.75 (0.06)  & 5.34 (0.01) & 35.84 (0.22) & 42.46 (0.05) \tabularnewline
		\multirow{1}{*}{HDEM }& 1.22 (0.01) & 8.60 (0.06) & 10.71 (0.06) & 1.93 (0.02) & 10.17 (0.12) & 15.34 (0.13) \tabularnewline
		\multirow{1}{*}{PEM }& 1.04 (0.01) & 6.67 (0.06) & 9.79 (0.04) & 1.26 (0.01) & 7.08 (0.05) & 10.91 (0.08) \tabularnewline
		\hline
		\multirow{1}{*}{\textbf{M2} } & \multicolumn{3}{c}{$p=400$} & \multicolumn{3}{c}{$p=1000$}\tabularnewline
		\hline 
		\multirow{1}{*}{Oracle }& 0.44 (0.003) & 12.29 (0.08) & 16.21 (0.02)& 0.45 (0.003) & 11.56 (0.07) & 15.97 (0.02)\tabularnewline 
		\multirow{1}{*}{PSEM }& 0.64 (0.003) & 13.67 (0.09) & 17.78 (0.04) & 0.70 (0.02) & 14.06 (0.09) & 18.15 (0.02) \tabularnewline               
		\multirow{1}{*}{GLLiM }& 6.74 (0.01) & 67.01 (0.20) & 33.95 (0.11) & 6.77 (0.01) & 68.60 (0.21) & 33.04 (0.10) \tabularnewline
		\multirow{1}{*}{Initial }& 2.79 (0.04) & 46.08 (0.24) & 44.03 (0.04)  & 2.99 (0.04) & 55.76 (0.24) & 44.75 (0.04) \tabularnewline
		\multirow{1}{*}{HDEM }& 1.26 (0.01) & 30.87 (0.12) & 22.58 (0.06)  & 1.96 (0.01) & 41.59 (0.22) & 29.65 (0.11)\tabularnewline
		\multirow{1}{*}{PEM }& 1.03 (0.004) & 23.63 (0.13) & 19.95 (0.04)  & 1.39 (0.01) & 33.36 (0.23) & 23.28 (0.09)\tabularnewline
		\hline
		\multirow{1}{*}{\textbf{M3} } & \multicolumn{3}{c}{$p=400$} & \multicolumn{3}{c}{$p=1000$}\tabularnewline
		\hline 
		\multirow{1}{*}{Oracle }& 0.56 (0.01) & 8.16 (0.06) & 10.09 (0.01) & 0.48 (0.01) & 7.60 (0.06) & 9.70 (0.02) \tabularnewline 
		\multirow{1}{*}{PSEM }& 0.67 (0.01) & 7.97 (0.06) & 11.78 (0.06) & 0.65 (0.01) & 7.81 (0.06) & 11.10 (0.05)\tabularnewline               
		\multirow{1}{*}{GLLiM }& 9.65 (0.01) & 68.41 (0.20) & 33.72 (0.11)  & 9.54 (0.01) & 67.52 (0.22) & 33.974 (0.10)\tabularnewline
		\multirow{1}{*}{Initial }& 5.02 (0.01) & 36.47 (0.20) & 43.41 (0.05)   & 5.30 (0.01) & 43.71 (0.22) & 43.79 (0.05) \tabularnewline
		\multirow{1}{*}{HDEM }& 1.54 (0.01) & 11.82 (0.07) & 14.12 (0.09)  & 2.07 (0.02) & 13.39 (0.10) & 17.68 (0.13) \tabularnewline
		\multirow{1}{*}{PEM }& 1.21 (0.01) & 10.67 (0.07) & 12.27 (0.06) & 1.39 (0.01) & 10.47 (0.07) & 12.10 (0.05) \tabularnewline
		\hline
		\multirow{1}{*}{\textbf{M4} } & \multicolumn{3}{c}{$p=400$} & \multicolumn{3}{c}{$p=1000$}\tabularnewline
		\hline 
		\multirow{1}{*}{Oracle }& 1.03 (0.03) & 5.14 (0.03) & 7.57 (0.01) & 1.35 (0.04) & 5.28 (0.05) & 7.53 (0.01) \tabularnewline 
		\multirow{1}{*}{PSEM }& 1.80 (0.05) & 5.21 (0.03) & 10.82 (0.12) & 3.56 (0.07) & 6.58 (0.06) & 16.53 (0.19)\tabularnewline               
		\multirow{1}{*}{GLLiM }& 17.38 (0.01) & 52.15 (0.27) & 33.79 (0.06) & 17.40 (0.01) & 54.10 (0.23) & 34.99 (0.06) \tabularnewline
		\multirow{1}{*}{Initial }& 13.16 (0.02) & 18.56 (0.10) & 49.62 (0.06) & 13.54 (0.01) & 18.48 (0.09) & 51.21 (0.06) \tabularnewline
		\multirow{1}{*}{HDEM }& 8.87 (0.09) & 8.93 (0.07) & 30.02 (0.25)  & 12.08 (0.07) & 19.28 (0.13) & 43.62 (0.24) \tabularnewline
		\multirow{1}{*}{PEM }& 2.43 (0.04) & 5.08 (0.03) & 11.49 (0.12) & 3.99 (0.06) & 6.32 (0.05) & 17.10 (0.19) \tabularnewline
		\hline
		
	\end{tabular}
}
\caption{\label{tab:binary} Average estimation errors based on 100 replicates (standard errors in parentheses). }
\end{table*}

\begin{figure*}[ht!]
\centering
\includegraphics[scale=0.56]{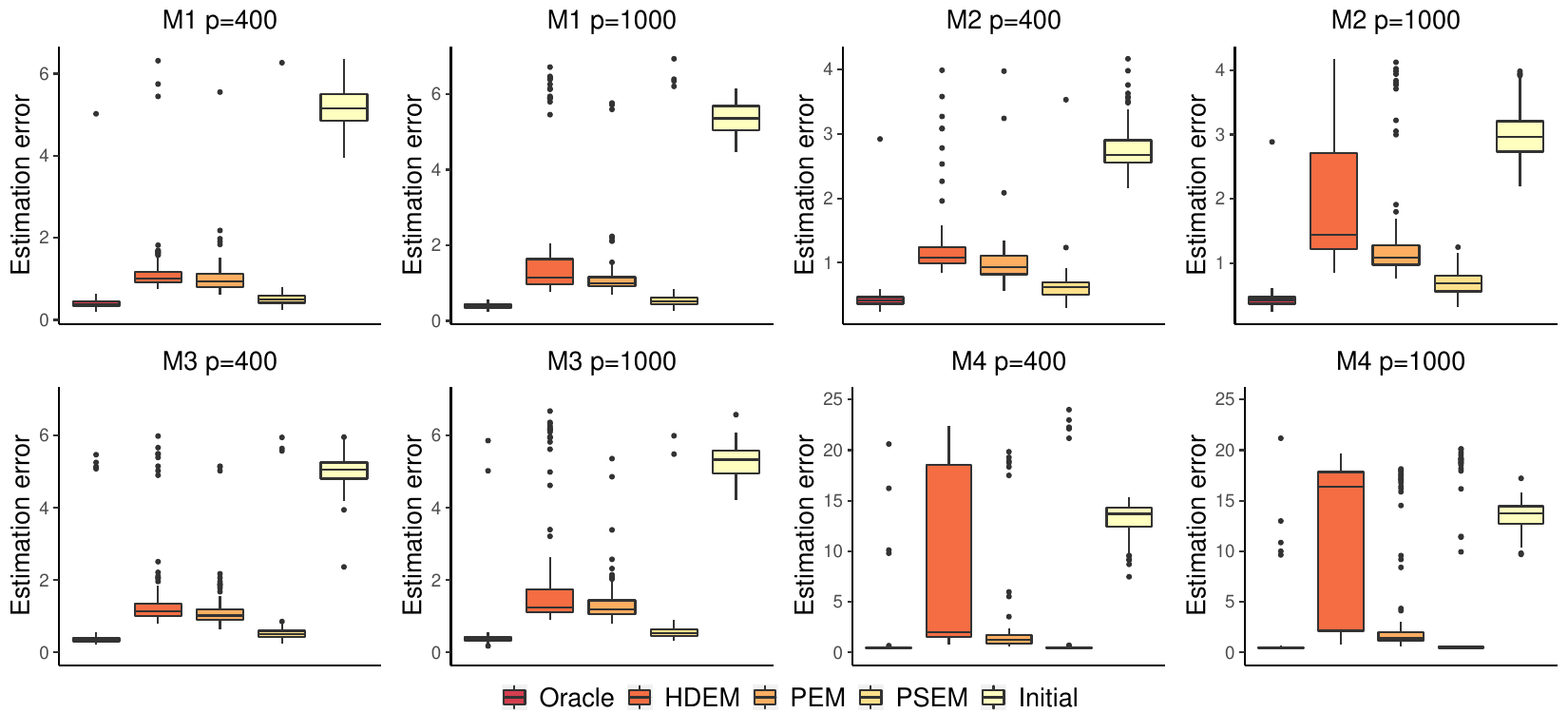}
\caption{The estimation error of $\bolbeta_k$'s based on 100 replicates. }\label{fig1}
\end{figure*}

\begin{table*}[ht!]
\centering
\resizebox{\textwidth}{!}{
	\begin{tabular}{ccccccccc}
		\hline 			
		\multirow{2}{*}{ }&\multicolumn{2}{c}{M1 $p=400$} & \multicolumn{2}{c}{M1 $p=1000$}& \multicolumn{2}{c}{M2 $p=400$} & \multicolumn{2}{c}{M2 $p=1000$}\tabularnewline&TPR & FPR &TPR & FPR &TPR & FPR &TPR & FPR \tabularnewline
		\hline 
		\multirow{1}{*}{Initial }& 83.2 (0.21) & 8.3 (0.05) & 76.4 (0.23) & 3.2 (0.02)& 88.3 (0.19) & 8.9 (0.05) & 79.6 (0.21) & 4.3 (0.03)\tabularnewline
		\multirow{1}{*}{HDEM }& 92.7 (0.09) & 1.4 (0.03) & 86.6 (0.15) & 1.5 (0.03)& 89.8 (0.01) & 1.9 (0.02)& 81.4 (0.25) & 1.2 (0.01)\tabularnewline
		\multirow{1}{*}{PEM }& 100 (0) & 0.9 (0.02)& 100 (0) & 0.7 (0.02)& 100 (0) & 1.2 (0.01) & 100 (0) & 0.8 (0.01)\tabularnewline
		\hline
		\multirow{1}{*}{ }&\multicolumn{2}{c}{M3 $p=400$} & \multicolumn{2}{c}{M3 $p=1000$}& \multicolumn{2}{c}{M4 $p=400$} & \multicolumn{2}{c}{M4 $p=1000$}\tabularnewline
		\hline 
		\multirow{1}{*}{Initial } & 89.1 (0.16) & 0.5 (0.05) & 73.4 (0.25) & 3.9 (0.03)& 47.5 (0.31) & 4.0 (0.05) & 38.3 (0.28) & 1.7 (0.02)\tabularnewline
		\multirow{1}{*}{HDEM } & 95.0 (0.06) & 2.6 (0.05)  & 80.2 (0.12) & 1.7 (0.03)& 83.4 (0.24) & 12.6 (0.15) & 64.7 (0.33) & 6.7 (0.06)\tabularnewline
		\multirow{1}{*}{PEM }& 100 (0) & 1.3 (0.02) & 100 (0) & 0.5 (0.01)& 100 (0) & 2.7 (0.11) & 100 (0) & 3.8 (0.09)\tabularnewline
		\hline
	\end{tabular}
}
\caption{\label{tab:binary1} Average true positive rate (TPR) and false positive rate (FPR) for variable selection based on 100 replicates (standard errors in parentheses). }
\end{table*}

The estimation results are summarized in Table \ref{tab:binary}, particularly for $\beta$'s estimation we further plot the results from 100 replicates in Figure \ref{fig1}. As expected, the \emph{oracle} method, i.e., the standard EM algorithm applied to the truly relevant $s=10$ predictors, has the best performance. On the other hand, the method that does not perform variable selection, GLLiM, failed for all the simulation settings. The proposed method, \emph{PEM}, has encouraging performances for all the simulation models; overall, it is comparable to the \emph{oracle} method and has a slight edge over the very recent method \emph{HDEM} \citep{zhang2020estimation}. The advantage of our method over \emph{HDEM}, which updates each $\beta_k$ separately via lasso penalized estimation in the M-step, can be explained by the group-wise penalization and joint estimation of all $\beta_k$'s in the M-step. In model (M4) with three mixtures, our approach of joint estimation has even more gains in estimation accuracy. 

We also note that Algorithm~\ref{alg:mlr} substantially improves over the \emph{Initial} values. Even when the initialization is quite far from the truth, the proposed algorithm improves over iterations in terms of all the evaluation criteria. This is partly explained by the guaranteed monotonicity of our iterative algorithm in the penalized log-likelihood (see Lemma S.1.). Finally, the post-selection EM algorithm, \emph{PSEM}, further removes the bias in the penalized estimation and improves our method one-step closer to the \emph{oracle} method. This is because our method has excellent variable selection accuracy, as summarized in Table \ref{tab:binary1}. 

\subsection{Simulations for misspecified $\sigma^2$ in Algorithm \ref{alg:mlr}}\label{addsimu}
In this section, we show some simulation results for using different $\sigma^2$ in Algorithm \ref{alg:mlr}. The data is generated from \eqref{mlr}.
We set $K=2$, $\omega_1=\omega_2=1/2$, $p=400$, $n=400$, and the true value of $\sigma^2$ to be 1. The nonzero coefficients of $\bolbeta_1$ is generated independently from $N(0,1)$; and the nonzero coefficients of $\bolbeta_2$ is $\bolbeta_{2j}=\bolbeta_{1j}+\delta\cdot\mathrm{sgn}(\bolbeta_{1j})$, $j=1,\ldots,p$. The signal strength $\delta$ takes value from $\{0.75, 1, 2\}$ and the input variance $\sigma^2$ in Algorithm \ref{alg:mlr} takes value from $\{0.5,0.75,1,1.5,2\}$.

\begin{figure}[ht!]
\begin{center}
	\includegraphics[scale=0.65]{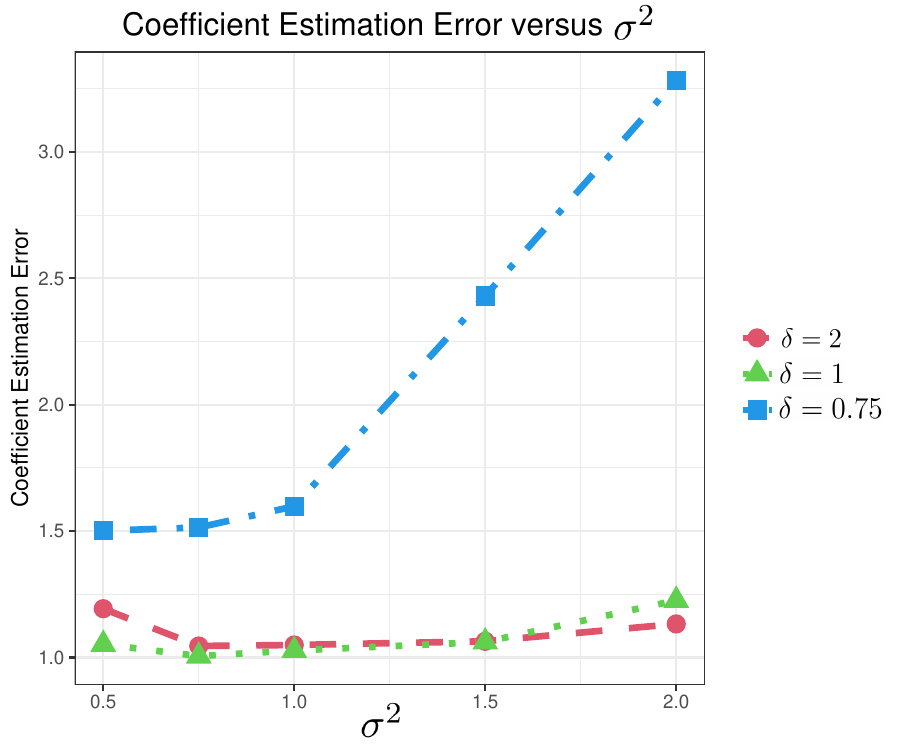}
\end{center}
\caption{Coefficient estimation errors for Algorithm \ref{alg:mlr} with different input $\sigma^2$ based on 100 replicates. }
\label{fig:sigma}
\end{figure}

\begin{table*}[ht!]
\centering
\renewcommand\arraystretch{1}
\begin{tabular}{cccccc}
	\hline 	
	& $\sigma^2=0.5$ & $\sigma^2=0.75$ & $\sigma^2=1$& $\sigma^2=1.5$ & $\sigma^2=2$ \tabularnewline 
	$\delta=2$  & 1.19 (0.08) & 1.05 (0.03) & 1.05 (0.03) & 1.06 (0.03) & 1.13 (0.07) \tabularnewline 
	$\delta=1$   & 1.05 (0.05) & 1.01 (0.04)  & 1.03 (0.04) & 1.06 (0.03) & 1.23 (0.05) \tabularnewline 
	$\delta=0.75$   & 1.50 (0.09) & 1.51 (0.10)  & 1.60 (0.09) & 2.43 (0.12) & 3.28 (0.07)\tabularnewline 
	\hline
\end{tabular}
\caption{\label{tab3} Coefficient estimation errors for Algorithm \ref{alg:mlr} with different input $\sigma^2$ based on 100 replicates. (standard errors in parentheses). }
\end{table*}

From Table \ref{tab3} and Figure \ref{fig:sigma}, when $\delta$ equals 1 and 2, using all the five $\sigma^2$ in Algorithm \ref{alg:mlr} returns good estimations. Although when 
$\sigma^2$ is 0.5 or 2, the performances are slightly worse, but can still capture the mixture information and give accurate enough estimates. The results are consistent with Theorem \ref{the:sigma}. For mixture linear regression models with a large signal-to-noise ratio, the choice of $\sigma^2$ has a minor influence on the performance of Algorithm \ref{alg:mlr}.
When $\delta$ is 0.75, the signal-to-noise ratio is not large enough, the performance of Algorithm \ref{alg:mlr} becomes worse, no matter the value of $\sigma^2$. This is mainly caused by the poor quality of the initialization. It deserves to notice that when $\sigma^2$ equals 0.5 or 0.75, the performance of the algorithm is slightly better than taking $\sigma^2=1$. This result indicates that when the signal-to-noise ratio is relatively small, Algorithm \ref{alg:mlr} can be less sensitive to the initialization if $\sigma^2$ takes a smaller value than the true one. 

\begin{algorithm}[ht!]
\caption{Group lasso penalized EM algorithm for model \eqref{mmlr}}
\label{alg:mmlr}
\begin{enumerate}
	\item[Input:] Initial values $\widehat{\omega}_k^{(0)}$, $\widehat{\beta}_k^{(0)}$, for $k=1,\cdots,K$, maximum iteration number $T$, data $\{X_i,Y_i; i=1,\ldots,n\}$, and initial tuning parameter	
	\begin{equation*}
		\begin{aligned}
			\lambda_n^{(0)}=C_1(\vert\widehat{ \omega}_1^{(0)}-\omega_1^*\vert\lor\Vert\widehat{\beta}_1^{(0)}-\beta_1^*\Vert_F\lor\cdots\lor\quad\Vert\widehat{\beta}_K^{(0)}-\beta_K^*\Vert_F)/{\sqrt{s}}+C_{\lambda}\sqrt{\logg(n)^2\logg(pq)/s}
		\end{aligned}
	\end{equation*}
	for some positive constants $C_1$ and $C_{\lambda}$. 
	
	\item[Iterate:] {For $t=0,\cdots,T-1$, do the following steps until convergence.} 
	\begin{itemize}
		
		\item For $i=1,\cdots,n$, let
		\begin{equation*}
			\begin{aligned}
				&\quad \widehat{\eta}_{i,k}(\widehat{\theta}^{(t)})={\widehat{\omega}_k^{(t)}}/\Big(\widehat{\omega}_k^{(t)}+\sum_{k'\neq k}\widehat{\omega}_{k'}^{(t)}\exp\big\{X_i^T(\widehat{\beta}_{k'}^{(t)}-\widehat{\beta}_{k}^{(t)})^T\bolSigma_y^{-1} \big(Y_i-{(\widehat{\beta}_k^{(t)}+\widehat{\beta}_{k'}^{(t)})^TX_i}/{2}\big)  \big\} \Big).
			\end{aligned}
		\end{equation*}
		
		\item For $k=1,\cdots, K$, update 
		\begin{equation*}
			\begin{aligned}
				&\widehat{\omega}_k^{(t+1)}=\frac{1}{n}\sum_{i=1}^n\widehat{\eta}_{i,k}(\widehat{\theta}^{(t)}),\\ &\widehat{\rho}_k^{(t+1)}=\frac{1}{n}\sum_{i=1}^n\big(\widehat{\eta}_{i,k}(\widehat{\theta}^{(t)})X_iY_i^T\bolSigma_y^{-1}\big),\\ &\widehat{\Sigma}_k^{(t+1)}=\frac{1}{n}\sum_{i=1}^n\big(\widehat{\eta}_{i,k}(\widehat{\theta}^{(t)})X_iX_i^T\big),
			\end{aligned}
		\end{equation*}
		and update $\widehat{\beta}_k^{(t+1)}$ by minimizing
		\begin{equation*}
			\begin{aligned}
				\sum_{k=1}^K\tr\big(\bolSigma_y^{-1}\beta_k^T\widehat{\Sigma}_k^{(t+1)}\beta_k\big)-2\sum_{k=1}^{K}\tr\big((\widehat{\rho}_k^{(t+1)})^T\beta_k\big)+\lambda_n^{(t+1)}\sum_{j=1}^p\sum_{l=1}^q\sqrt{\sum_{k=1}^K(\beta_{k})_{jl}^2}.
			\end{aligned}
		\end{equation*}
		with $\lambda_n^{(t+1)}=\kappa\lambda_n^{(t)}+C_{\lambda}\sqrt{\logg(pq)\logg(n)^2/n}$, where $\kappa\in(0,1/2)$ and $(\beta_{k})_{jl}$ is the $(j,l)$-th element of $\bolbeta_k$.
	\end{itemize}
	\item[Output:] $\widehat\theta^{(t+1)}=\{\widehat\omega_1^{(t+1)},\ldots,\widehat\omega_K^{(t+1)},\widehat\beta_1^{(t+1)},\ldots,\widehat\beta_K^{(t+1)}\}$
\end{enumerate}
\end{algorithm}

\section{Extension to multivariate mixture linear regression}\label{MMLR}
\subsection{Model and algorithm}

In this section, we consider a generalization of the mixture linear regression model, where $Y$ is a $q$-dimensional response. Specifically, we consider the model
\begin{equation}\label{mmlr}
\mbbP(W=k)=\omega_k,\quad Y\mid (X,W=k)\sim N(\beta_k^T X,\Sigma_y),
\end{equation}
where $\beta\in\mbbR^{p\times q}$ and $\Sigma_y\in \mbbR^{q\times q}$ is a symmetric and positive definite matrix. We will answer two questions: \textit{Does the penalized EM algorithm still work and have similar convergence results for multiple response cases? If so, what is the advantage of considering $q$ responses simultaneously than $q$ mixture linear regression problems separately?}  

As the mixture linear regression model, to handle the high-dimensionality of $X$, we assume the matrix coefficients $\beta_k$, $k=1,\cdots,K$, to be sparse, and the sparsity patterns are the same for all the mixtures. We redefine $S=\{(i,j): (\beta_{k})_{ij}\neq 0,\ i=1\cdots,p,\ j=1,\cdots,q  \}$ and $s=\vert S\vert _0$ be the cardinality of $S$. We allow $p\gg n$ and $q$ to grow linearly with $s$.

In model \eqref{mmlr}, our parameter of interest is $\boltheta=\{\omega_1,\cdots,\omega_k,\bolbeta_1,\cdots,\bolbeta_K\}$. The error covariance $\Sigma_y$ is treated as a known parameter. The penalized EM algorithm for solving \eqref{mmlr} is analogous to Algorithm \ref{alg:mlr} and is summarized in Algorithm \ref{alg:mmlr}. Similar to Algorithm \ref{alg:mlr}, the formula of $\lambda_n^{(t)}$ is mainly for theoretical consideration. In practice, we fix $\lambda_n^{(t)}=\lambda$ over all iterations.

\subsection{Theoretical analysis}
In this section, we consider the theoretical studies for Algorithm \ref{alg:mmlr}. We assume that $\bolSigma_y$ is a known parameter and without loss of generality, we assume $\bolSigma_y=\mbI_q$ (If $\bolSigma_y$ is not identity, treating $\bolSigma_y^{-1/2}\bolbeta_k$ as the target parameter will turn the problem to the case $\bolSigma_y=\mbI_q$). Recall that, even for the univariate-response mixture linear regression, theoretical analysis for unknown variance $\sigma^2$ is extremely challenging and has not been considered in high dimensions. As a result, for multiple-response mixture linear regression, it makes sense to assume that $\bolSigma_y$ is known as a necessary simplification. As such, we focus on the coefficients $\bolbeta_k$ as parameters of interest. To the best of our knowledge, there is little work for multiple-response mixture regression even in this simplified context. 


We first define the following additional notations.
For a matrix $\mbA\in\mbbR^{p\times q}$, $\Vert\mbA\Vert_F$ is its Frobenius norm and $\Vert\mbA\Vert_{F,s}=\sup_{\bolmu\in\mbbR^{p\times q},\Vert\bolmu\Vert_F=1,\vecc(\bolmu)\in\bolGamma(s)}\langle\mbA,\bolmu\rangle_F$. 
The true parameter space we consider now is
\begin{equation*}
\bolTheta^*=\{\boltheta^*: \omega_1^*\in (c_w,1-c_w), \Vert\vecc(\bolbeta_k^*)\Vert_0\leq s, \Vert\bolbeta_k^*\Vert_F\leq M_b,\mathrm{for}\ k=1,2 \}, 
\end{equation*}
and the signal-to-noise ratio is re-defined as $\bolDelta=\sqrt{\tr((\bolbeta_2^*-\bolbeta_1^*)^T\bolSigma (\bolbeta_2^*-\bolbeta_1^*) )}/q$. Then we define $d_{F,s}\big( M(\boltheta_1),M(\boltheta_2)\big)$ and $d_{F}\big( M(\boltheta_1),M(\boltheta_2)\big)$ as
\begin{equation*}
\begin{aligned}
	&\max_{k=1,2}\{\vert\omega_k(\boltheta_1)-\omega_k(\boltheta_2)\vert \lor\Vert\bolrho_k(\boltheta_1)-\bolrho_k(\boltheta_2)\Vert_{F,s}\lor\Vert(\bolSigma_k(\boltheta_1)-\bolSigma_k(\boltheta_2))\bolbeta_2^*\Vert_{F,s}\},\\
	&\max_{k=1,2}\{\vert\omega_k(\boltheta_1)-\omega_k(\boltheta_2)\vert \lor\Vert\bolrho_k(\boltheta_1)-\bolrho_k(\boltheta_2)\Vert_{F}\lor\Vert(\bolSigma_k(\boltheta_1)-\bolSigma_k(\boltheta_2))\bolbeta_2^*\Vert_{F}\},
\end{aligned}
\end{equation*}
respectively, and the new constriction basin $\calB_{con}(\boltheta)$ as
\begin{equation*}
\calB_{con}(\boltheta)=\{\boltheta:\omega_k\in (c_0,1-c_0), \Vert\bolbeta_k-\bolbeta_k^*\Vert_F\leq C_b\Delta,\  \vecc(\bolbeta_k-\bolbeta_k^*)\in\Gamma(s),  \ \mathrm{for}\ k=1,2\}.
\end{equation*}

We require the same technical conditions as the mixture linear regression with a modification to Condition (C3):
\begin{enumerate}[(C3')]
\item The new signal-to-noise ratio $\Delta=\sqrt{\tr((\bolbeta_2^*-\bolbeta_1^*)^T\bolSigma (\bolbeta_2^*-\bolbeta_1^*) )}/q>C_1(c_0)$ for a constant $C_1(c_0)$ only depends on $c_0$, and $C_b<C_2(c_0,M_2)$ for a constant $C_2(c_0,M_2)$ only depends on $c_0,M_2$.
\end{enumerate}
The new signal-to-noise ratio reveals a fundamental difference of considering $q$ responses together than $q$ mixture linear regressions separately. To see this, define $\beta_{k,j}$ as the $j$-th column of $\beta_k$, i.e, the regression coefficient of $(Y_j,\mbX)$. Further define the individual signal-to-noise-ratio as $\Delta_j=\sqrt{(\bolbeta_{1,j}-\bolbeta_{2,j})^T\bolSigma (\bolbeta_{1,j}-\bolbeta_{2,j})}$. If we apply Algorithm~\ref{alg:mlr} to $(Y_j,X)$ for $j=1,\ldots,q$, then we need all of $\Delta_j$ to be bounded below, as discussed in Section \ref{sec:the}. In other words, when some $\Delta_j$ are small and the two mixtures are not well separated on this coordinate, Algorithm~\ref{alg:mlr} may not be able to estimate the corresponding coefficients. However, if we consider all the responses simultaneously in multivariate mixture linear regression, we only require $\Delta=\sqrt{\sum_{j=1}^q \Delta_j^2}/q$, the average of $q$ signal-to-noise ratios, to be greater than a constant. We can allow some responses to have small signal-to-noise ratios without rendering Algorithm~\ref{alg:mmlr} inapplicable. This is because all the response shares the common latent structure. The responses with  large signal-to-noise ratios can help enhance the identification of the mixtures for the responses with small signal-to-noise ratios. 

As a brief numerical illustration, we generated simulated datasets from \eqref{mmlr} with $K=2$, $n=400$, $p=100$, $s=10$, $q=2$, $\bolSigma_y=\mbI_q$ and $[\bolSigma]_{ij}=0.3^{\vert i-j\vert}$. The first 5 rows of $\bolbeta_k$, $k=1,2$, are set to be nonzero. More specifically, the nonzero coefficients of the first and second columns of $\bolbeta_1$ are generated from $N(0,1)$; and the nonzero coefficients of $\bolbeta_2$ are $(\bolbeta_2)_{j1}=(\bolbeta_1)_{j1}+2\cdot\mathrm{sgn}((\bolbeta_1)_{j1})$ and $(\bolbeta_2)_{j2}=(\bolbeta_1)_{j2}+\delta\cdot\mathrm{sgn}((\bolbeta_1)_{j2})$. The signal strength is 2 for the first response and $\delta$ for the second one. We vary $\delta$ from 0.5 to 2 with increments of 0.5. 
\begin{figure*}[ht!]
\centering
\includegraphics[scale=0.451]{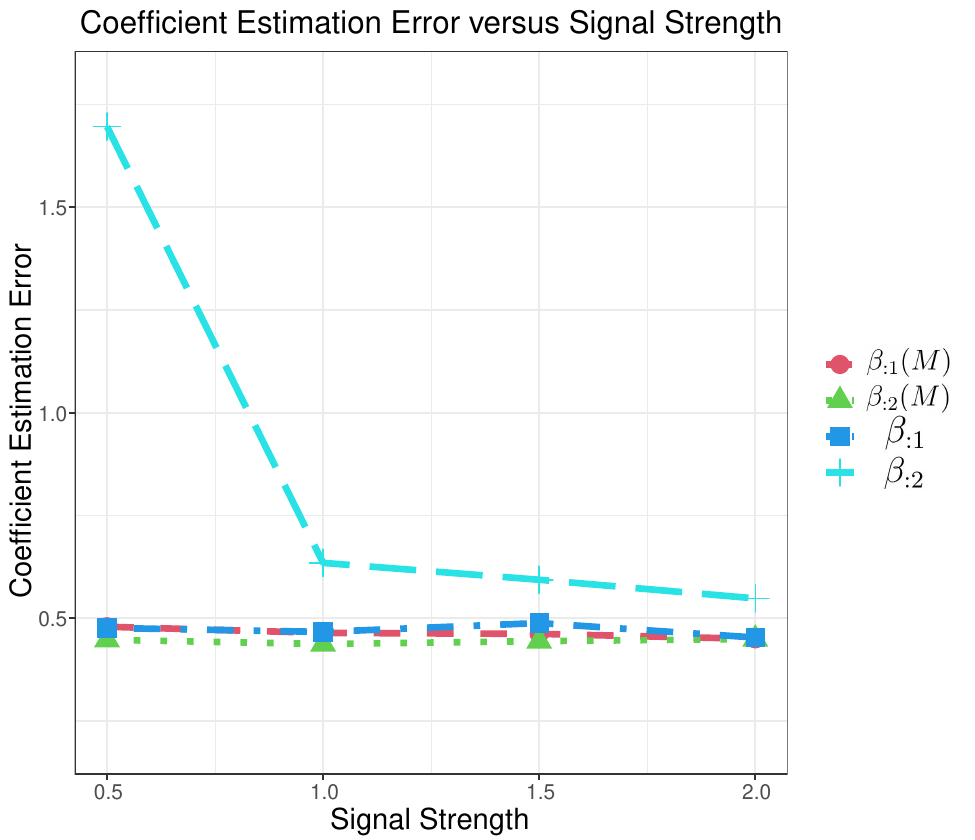}
\includegraphics[scale=0.46]{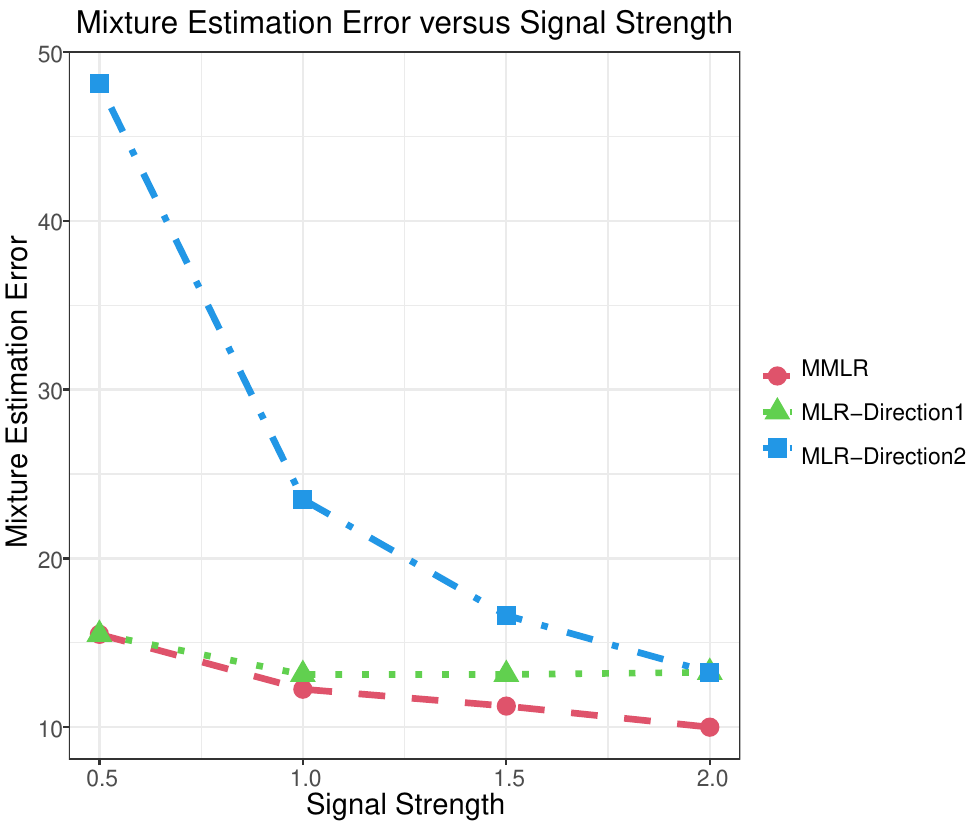}
\caption{The reported results for each signal strength are medians of 100 replicates. In the left panel, $\bolbeta_{:j}(M)$ represents the estimation error for the $j$-th column of $\beta$ using Algorithm \ref{alg:mmlr}; $\bolbeta_{:j}$ represents the estimation error for the $j$-th column of $\beta$ using Algorithm \ref{alg:mlr} separately for each response. In the right panel, MMLR represents the estimation error for the mixtures using Algorithm \ref{alg:mmlr}; MLR-Response1 and MLR-Response2 represent the estimation error for the mixtures using Algorithm \ref{alg:mlr} separately for the first and second element of the response, respectively.}\label{fig2}
\end{figure*}
From Figure \ref{fig2}, we can see that the coefficient estimation errors for Algorithm \ref{alg:mmlr} remain at a low level for all signal strengths and the mixture estimation error gains slightly improvement when the signal strength of the second response increases. When the signal strength of the second response is weak, the strong signal of the first response assists the recovery of the mixtures and helps to increase the coefficient estimation accuracy of the second coefficient substantially. As a comparison, when using Algorithm \ref{alg:mlr} separately for each response, the estimation errors for the second coefficient and mixture increase quickly when the signal strength decreases. In terms of the mixture estimation error, it even performs like random guessing when the signal strength is 0.5. 

The following lemmas and theorem rigorously show the statistical convergence for Algorithm \ref{alg:mmlr}.
\begin{lemma}\label{contraction2}
Under conditions (C1),(C3'), if $q=O(s)$ and $\theta\in\calB_{con}(\boltheta^*)$, then
\begin{equation*}
	d_{F}\big(M(\boltheta),M(\boltheta^*)\big)\leq \kappa_0\big( \vert\omega_1(\boltheta)-\omega_1^*\vert\lor\Vert\bolbeta_1-\bolbeta_1^*\Vert_F\lor\Vert\bolbeta_2-\bolbeta_2^*\Vert_F \big).
\end{equation*}
for some $0<\kappa_0<\frac{1}{2\lor (64/\tau_0)}$.
\end{lemma}
\begin{lemma}\label{concentation2}
Suppose that $\boltheta^*\in\bolTheta^*$. Under condition (C1), there exists a constant $C_{con}>0$, such that with probability at least $1-4(pq)^{-1}$,
\begin{equation*}
	\begin{aligned}
		&\sup_{\boltheta\in \calB_{con}(\boltheta^*)}d_{F,s}(M(\boltheta),M_n(\boltheta))\leq C_{con} \sqrt{\frac{s\logg(n)^2\logg(pq)}{n}}
	\end{aligned}
\end{equation*}
\end{lemma}
Lemmas \ref{contraction2} \& \ref{concentation2} can be interpreted similarly to Lemmas 1 and 2, respectively. However, since we are now interested in multivariate response, we further assume condition $q=O(s)$ in addition to the technical conditions (C1) to (C4). The assumption on $q$ results from the signal-to-noise ratio. Note that $\Delta=O(s/q)$. In order for the condition $\Delta>C_1$ to hold, we need $q=O(s)$. However, when this assumption holds, $q$, i.e, the dimensionality of $\mbY$, has a relatively small effect on the estimation error, as we only have an additional factor of $\logg{q}$ in Lemma~\ref{concentation2} in comparison to Lemma~\ref{concentation}.

Combine Lemmas \ref{contraction2} \& \ref{concentation2} and we have the following theorem.

\begin{theorem}\label{the:2}
Under conditions (C1), (C2), (C3'), and(C4) and $q=O(s)$, there exists a constant $0<\kappa<1/2$, such that $\hatbolbeta_k^{(t+1)}$ obtained by Algorithm \ref{alg:mmlr} satisfies, with probability greater than $1-4(pq)^{-1}$,
\begin{equation*}
	\begin{aligned}
		\Vert\hatbolbeta_k^{(t+1)}-\bolbeta_k^*\Vert_F&=O\Big(\kappa^t(\vert\widehat{ \omega}_1^{(0)}-\omega_1^*\vert\lor\Vert\hatbolbeta_1^{(0)}-\bolbeta_1^*\Vert_F\lor\Vert\hatbolbeta_2^{(0)}-\bolbeta_2^*\Vert_F) +\sqrt{\frac{s\logg(n)^2\logg(pq)}{n}}\Big).
	\end{aligned}
\end{equation*}
Consequently, for $t\geq \{-\logg(\kappa)\}^{-1}\logg\{n(\vert\widehat{ \omega}_1^{(0)}-\omega_1^*\vert\lor\Vert\hatbolbeta_1^{(0)}-\bolbeta_1^*\Vert_F\lor\Vert\hatbolbeta_2^{(0)}-\bolbeta_2^*\Vert_F)\}$,
\begin{equation*}
	\Vert\hatbolbeta_k^{(t+1)}-\bolbeta_k^*\Vert_F=O\Big(\sqrt{\frac{s\log(n)^2\logg(pq)}{n}}\Big).
\end{equation*}
\end{theorem}
Again, the convergence rate for multivariate mixture linear regression model is similar to that for the mixture linear regression model with an additional $\logg(q)$ term. The convergence rate grows slowly as a function of $q$. However, theoretical requirement $\Delta>C_1$ and practical consideration about the estimation for $\Sigma_y$ restrict $q$ to grow linearly with $s$ and slower than $n$. It would be interesting to develop a method that allows $q$ to grow faster, which we leave a challenging future topic. 

\section{Real data analysis} \label{sec:real}
The cancer cell line encyclopedia (CCLE) dataset contains 8-point dose-response curves for 24 chemical compounds across over 400 cell lines, which is a publicly available dataset at \url{ www.broadinstitute.org/ccle} and is also considered in \cite{li2019drug}. Because the cell lines are not consistent for different chemical compounds, we consider chemical compounds: Lapatinib, AZD6244, and PD-0325901, three popular chemical compounds for cancer treatment. Analogous to \cite{li2019drug}, we use the area under the dose-response curve, also known as the activity area, to measure the sensitivity of a drug for each cell line. Besides the drug information, the data also contains the expression data of 18,926 genes for each cell line. Aiming to identify the genes sensitive to the chemical compounds, we treat the active area as the response and the gene expressions as the predictor. After identifying cell lines treated using all three chemical compounds, we get the sample size $n=490$. We keep $p=500$ genes that are highly correlated with three responses (for each gene, we calculated the sum of its absolute correlations with the three responses). Due to the complexity of cancer, we expect the data to be heterogeneous.

We consider $K=2, 3, 4, 5, 6$, and $8$ for the proposed algorithms and denote individual mixture linear regression using Algorithm \ref{alg:mlr} as MLR and multivariate mixture linear regression using Algorithm \ref{alg:mmlr} as MMLR. In addition, we consider the methods used in the simulation studies including HDEM, Initial method, and PSEM and LASSO regression as competitors. We repeatedly split the whole data set into $1:4$ ratios, and use $1/5$ of the data as the testing samples and the rest as the training samples. The process is repeated 100 times. We report the mean squared prediction errors. Please see Table \ref{tab4} for the overall prediction errors and Figures \ref{cclebox} and \ref{cclebox1} for the boxplots. More detailed comparison results for each response are presented in Section S.8 of the Supplement.

When $K=5$, the overall mean squared prediction errors for MMLR is the smallest. However, from Figure \ref{cclebox}, we can see that MMLR is insensitive to $K$. When $K=4,5$ and $6$, the mean squared prediction errors are very close. As a comparison, MLR is more sensitive to $K$ and achieves the best prediction error when $K=4$. Two possible reasons make MMLR perform better than MLR. The first reason is that considering the three responses together results in a larger signal-to-noise ratio, which facilitates the mixture identification and coefficients estimation. The second reason is that, unlike MMLR, the estimated mixtures are different for the three responses in MLR, which may reduce the estimation accuracy. Although MLR performs slightly worse than MMLR, it still outperforms the other competitors, especially the LASSO regression. It implies that the data set is quite heterogeneous, the proposed mixture linear regression approach improves the performance of a single linear regression significantly. 

When $K=4$ and $6$, we also consider the estimated label for each observation based on the full data set. The number of observations in each estimated mixture is reported in Table \ref{tab5}. Based on this table, we recommend using $K=2$, 3 and 3 for MLR to the three responses, respectively, and using $K=5$ for MMLR. Except from the first response, the recommended $K$'s make MLR and MMLR achieve the lowest or among the lowest prediction error. For the first response, MLR achieves the smallest prediction error when $K=6$. However, it only separates the observations into 2 mixtures. Note that the estimated label for the $i$th observation is $\argmax_{k=1,\cdots,K}\widehat{\eta}_{i,k}(\widehat{\theta})$, it can happen that no observation is classified into some mixtures. For the first response, the two mixtures are more balanced with the increase of $K$, which may be one reason makes the prediction error smaller. Another reason is that the relationship between the first response and the predictors may be complicated and non-linear, which can not be captured by a mixture linear regression model. Hence, MLR tends to use more mixtures to approximate the true relationship. Another phenomena we observe is that MMLR uses more mixtures than MLR. This happens since the true mixtures can be different for the three responses. By separating the samples into more mixtures, MMLR becomes more accurate and robust.

\begin{table*}[ht!]
	\centering
	\renewcommand\arraystretch{1}
	\resizebox{15.2cm}{1.45cm}{
		\begin{tabular}{c|cccccc}
			\hline
			K &MLR  &MMLR&HDEM&Initial&PSEM &LASSO \tabularnewline 
			2&0.413 (0.005) & 0.388 (0.004) & 0.420 (0.004) & 0.612 (0.009) &1.461 (0.020)    &\multirow{5}{*}{0.948 (0.026)} 
			\tabularnewline 
			3&0.275 (0.004) &0.276 (0.004)&0.347 (0.004) &0.493 (0.010)&1.376 (0.021)
			\tabularnewline 
			4&0.254 (0.005) & 0.238 (0.004) & 0.324 (0.004) & 0.443 (0.008) &1.079 (0.020)
			\tabularnewline 
			5&0.276 (0.007) &\textbf{0.237} (0.004)&0.315 (0.005) &0.416 (0.008)&0.569 (0.012)
			\tabularnewline 
			6&0.302 (0.007) & 0.239 (0.004) & 0.310 (0.004) & 0.395 (0.007) &0.412 (0.008)
			\tabularnewline 
			8&0.347 (0.006) & 0.243 (0.004) & 0.310 (0.005) & 0.357 (0.006) &0.375 (0.007)
			\tabularnewline 
			\hline
		\end{tabular}
	}
	\caption{\label{tab4} Mean squared prediction errors based on 100 replicates (Standard errors are given in parenthesis). }
\end{table*}

%

\begin{figure*}[ht!]
	\centering
	\includegraphics[scale=0.3]{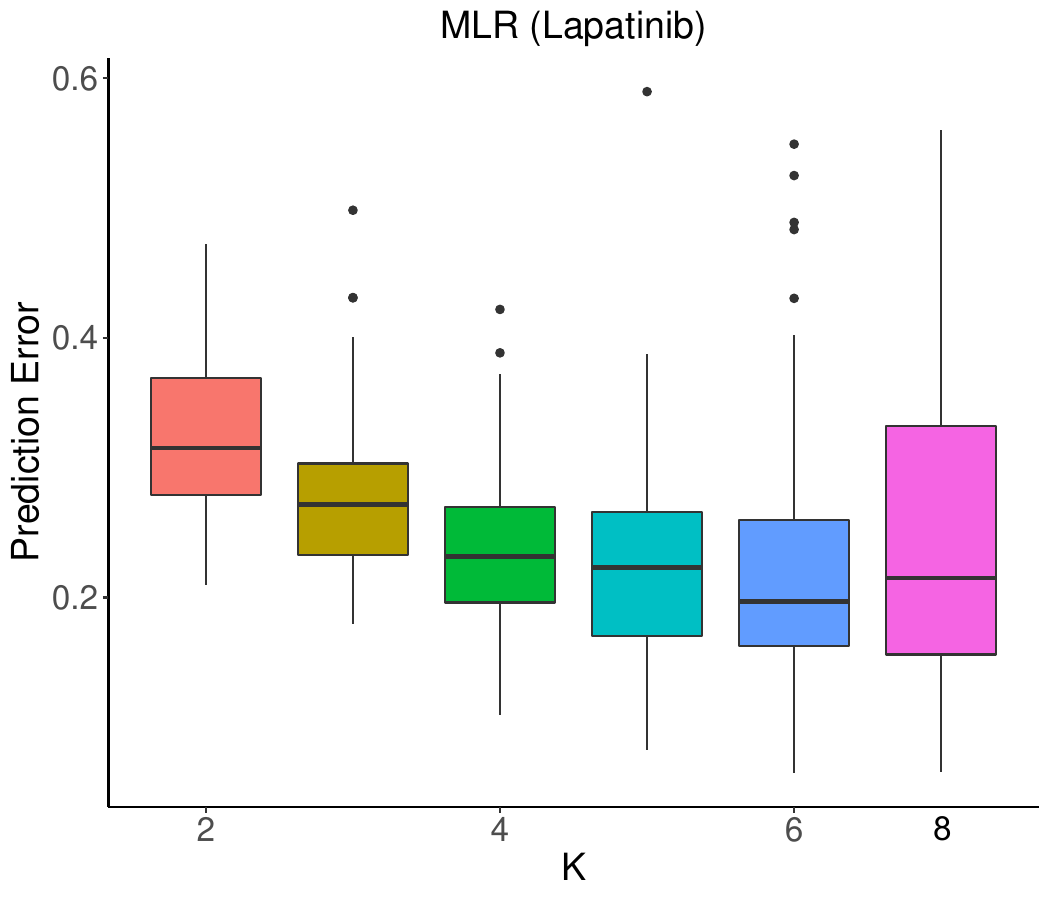}
	\includegraphics[scale=0.3]{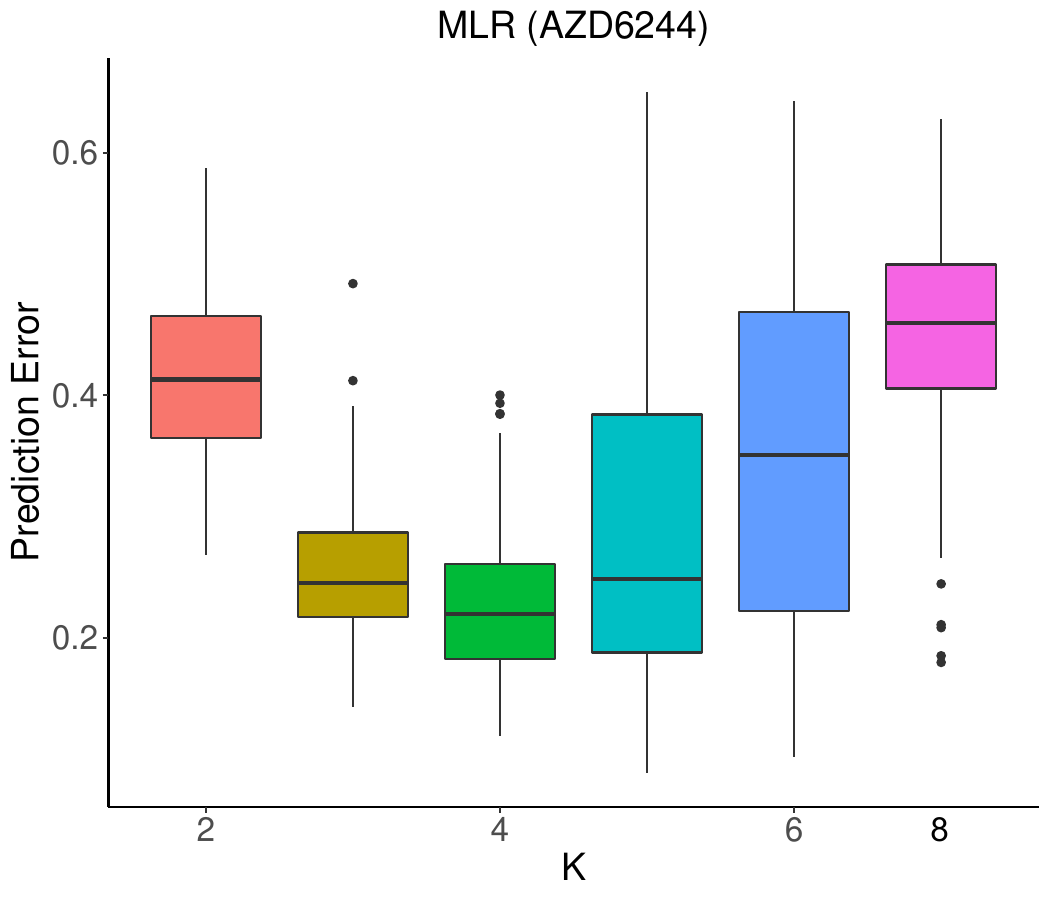}
	\includegraphics[scale=0.3]{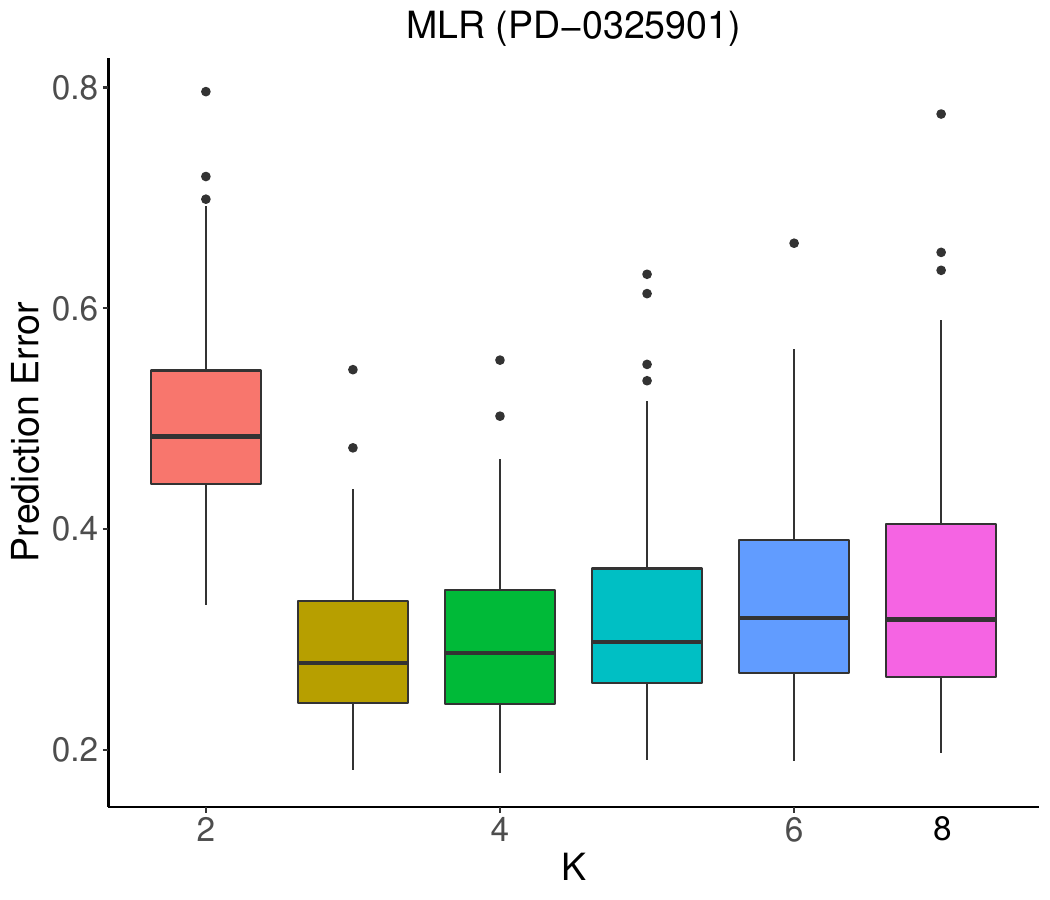}
	\includegraphics[scale=0.27]{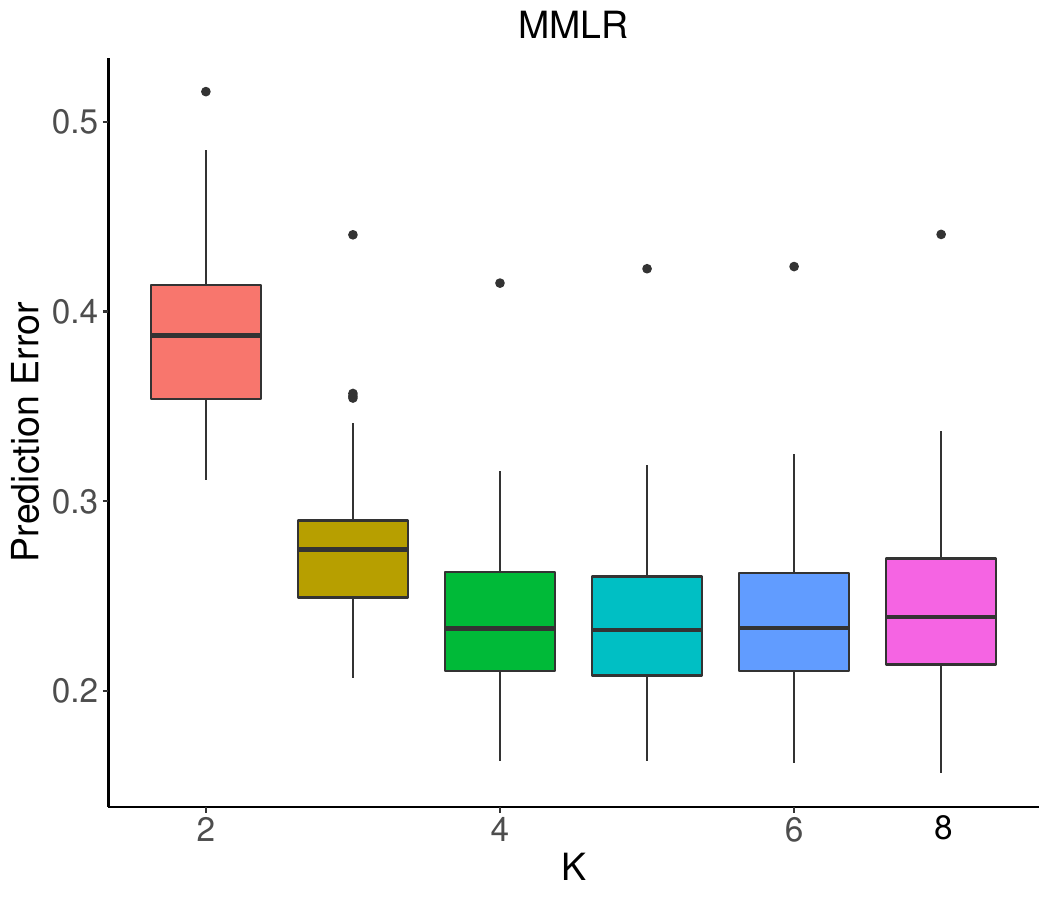}
	\caption{CCLE data: Boxplot for the mean squared prediction error }\label{cclebox}
\end{figure*}

\begin{figure*}[ht!]
	\centering
	\includegraphics[scale=0.36]{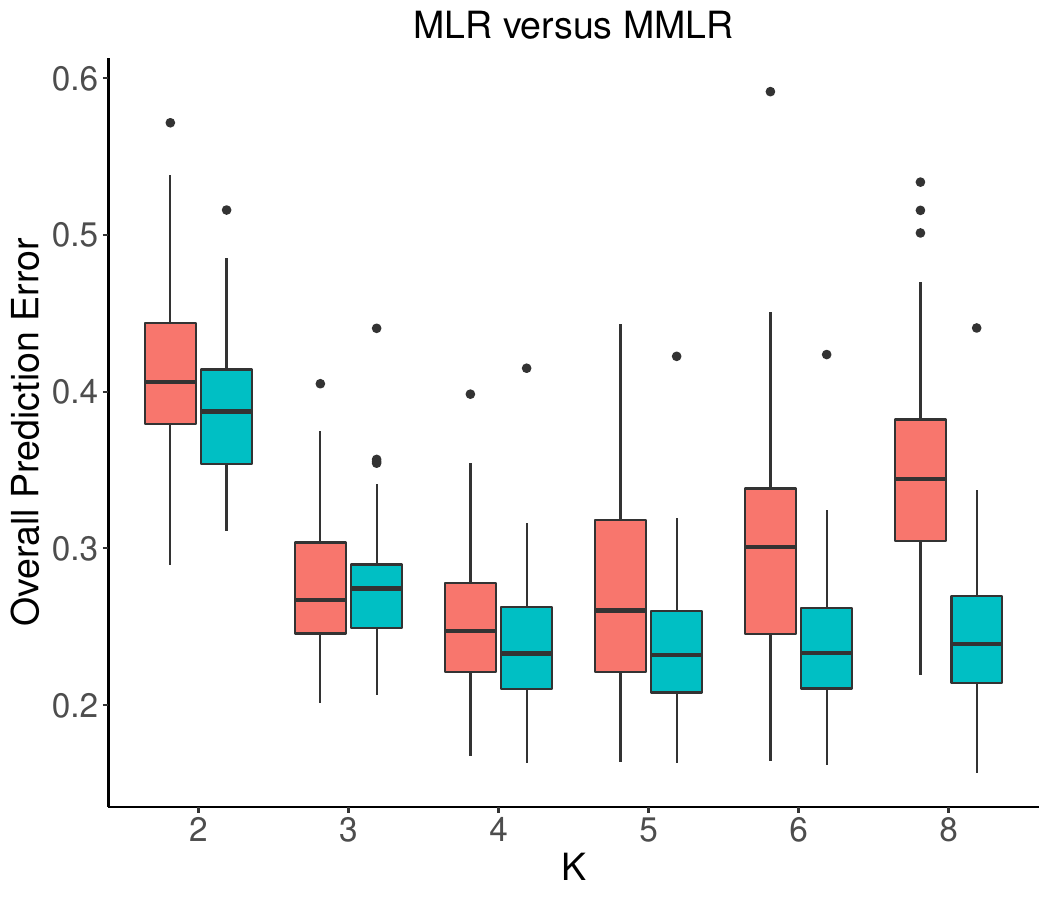}
	\caption{Overall prediction error for MLR and MMLR. The red boxes represent results for MLR.}\label{cclebox1}
\end{figure*}

\begin{table*}[ht!]
	\centering
	\renewcommand\arraystretch{1}
	\begin{tabular}{c|cccc|cccccc}
		\hline
		& \multicolumn{4}{c|}{K=4} & \multicolumn{6}{c}{K=6} \tabularnewline 
		\multirow{1}{*}{MLR (Lapatinib)}  &0 &0 &147 &343&0&0&0&0&162&328 \tabularnewline
		\multirow{1}{*}{MLR (AZD6244)}   &0 &87&162 &241 &0&0&0&85&187&218 \tabularnewline 
		\multirow{1}{*}{MLR (PD-0325901)}   &0 &105 &170 &215&0&0&0&101&180&209\tabularnewline 
		\multirow{1}{*}{MMLR}   &76&128&138 &148&0&48&73&104&114&151\tabularnewline 
		\hline 
	\end{tabular}
	\caption{\label{tab5} Numbers of samples in each estimated mixtures.}
\end{table*}


\section{Discussion}\label{sec:dis}
The paper studies a group lasso penalized EM algorithm for high-dimensional mixture linear regression. We obtained an encouraging non-asymptotic convergence rate without data splitting and under a general model setting. Although our theory is for normally distributed predictors and two-mixtures regression, the penalized EM algorithm is applicable to essentially any predictors and to more than two mixtures. We then extended the mixture linear regression model and the penalized EM algorithm to the multivariate response case and established its non-asymptotic theory.

The theoretical study is currently for the model with two mixtures. It provides insights and foundations to the theoretical studies for cases with $K\geq 3$. We leave it as a future work. Besides, the error is assumed to follow a normal distribution, it is interesting to considering more general error assumptions, such as the student t distribution, to achieve robustness. Finally, it is also a meaningful future work to extend the studies to generalized mixture linear regression models.

\baselineskip=15pt
\bibliography{ref_tensor1}
\bibliographystyle{agsm}

\end{document}